\documentstyle[amsmath,amssymb,11pt%
]{article}
\renewcommand{\title}[1]{\leftline{\Large\bf #1}\par\medskip}
\renewcommand{\author}[1]{\medskip{\large #1}\par\medskip}

\newcommand{\om}{\omega}

\newcommand{\rom}[1]{{\rm #1}}

\allowdisplaybreaks

\makeatletter\@addtoreset{equation}{section}\makeatother

\begin{document}

\setcounter{page}{1} \setcounter{section}{0} \thispagestyle{empty}

\newtheorem{definition}{Definition}[section]
\newtheorem{remark}{Remark}[section]
\newtheorem{proposition}{Proposition}[section]
\newtheorem{theorem}{Theorem}[section]
\newtheorem{corollary}{Corollary}[section]
\newtheorem{lemma}{Lemma}[section]
\newtheorem{example}{Example}

\newcommand{\skl}{\overset{(k)}{\diamondsuit}}
\newcommand{\D}{{\cal D}}
\newcommand{\N}{{\Bbb N}}
\newcommand{\C}{{\Bbb C}}
\newcommand{\Z}{{\Bbb Z}}
\newcommand{\R}{{\Bbb R}}
\newcommand{\Rp}{{\R_+}}
\newcommand{\eps}{\varepsilon}

\newcommand{\AS}{\operatorname{AS}}

\newcommand{\fii}{\varphi}
\newcommand{\Ker}{\operatorname{Ker}}
\newcommand{\supp}{\operatorname{supp}}
\newcommand{\la}{\langle}
\newcommand{\ra}{\rangle}
\newcommand{\const}{\operatorname{const}}
\newcommand{\ddGamma}{\overset{{.}{.}}{\Gamma}_X}

\newcommand{\ho}{\widehat\otimes}
\newcommand{\ot}{\otimes}

\renewcommand{\emptyset}{\varnothing}
\newcommand{\Formsg}{{\cal F}\Omega^{n+1}}
\newcommand{\Func}{{\cal FC}}
\newcommand{\FC}{{\cal F}C_{\mathrm b}^\infty({\cal D},\Gamma_X)}

\newcommand{\di}{\partial}
\renewcommand{\div}{\operatorname{div}}

\begin{center}{\Large \bf
Laplace operators in deRham complexes }\\[2mm] {\Large\bf
associated with measures}
\\[3mm] {\Large\bf on configuration spaces}\\[4mm]  \large Sergio
Albeverio,  Alexei Daletskii,
  Yuri Kondratiev,  Eugene Lytvynov

\end{center}

\begin{abstract}
Let $\Gamma _X$ denote the space of all locally finite
configurations in a
complete, stochastically complete, connected, oriented Riemannian manifold $%
X $, whose volume measure $m$ is infinite. In this paper, we
construct and study spaces $L^2_\mu\Omega^n$ of differential
$n$-forms over $\Gamma_X$ that are square integrable with respect
to a probability measure $\mu$ on $\Gamma_X$. The measure $\mu$ is
supposed to satisfy the  condition $\Sigma_m'$ (generalized Mecke
identity) well known in the theory of point processes. On
$L^2_\mu\Omega^n$, we introduce bilinear forms  of Bochner and
deRham type. We prove their closabilty and call the generators of
the corresponding closures  the Bochner and deRham Laplacian,
respectively. We prove that both operators contain in their domain
the set of all smooth local forms. We show that, under a rather
general assumption on the measure $\mu$, the space of all
Bochner-harmonic $\mu$-square integrable forms on $\Gamma_X$
consists only of the zero form. Finally, a Weitzenb\"ock type
formula connecting the Bochner and deRham Laplacians is obtained.
As examples, we consider (mixed) Poisson measures, Ruelle type
measures on $\Gamma_{\R^d}$, and Gibbs measures in the low
activity--high temperature regime, as well as Gibbs measures with
a positive interaction potential on $\Gamma_X$.
\end{abstract}


\thispagestyle{empty} \setcounter{page}{0}\newpage

 \tableofcontents

\section{Introduction}

Let $\Gamma _X$ denote the space of all locally finite
configurations in a
complete, stochastically complete, connected, oriented Riemannian manifold $%
X $ of infinite volume. The growing interest in geometry and
analysis on the configuration spaces $\Gamma _X$ can be explained
by the fact that these naturally appear in different problems of
statistical mechanics, quantum physics, and the theory of point
processes. In \cite{AKR-1,AKR0,AKR1}, an approach to the
configuration spaces as infinite-dimensional manifolds was
initiated. This approach was motivated by the theory of
representations of diffeomorphism groups, see \cite{GGPS, I,VGG}
(these references as well as \cite{AKR1, AKR3} also contain
discussion of relations with quantum physics). We refer the reader
to \cite{AKR2, AKR3,Lipsche, Ro} and references therein for
further discussion of analysis on the configuration spaces and
applications.

On the other hand, stochastic differential geometry of
infinite-dimensional manifolds, in particular, their (stochastic)
cohomologies and related questions (Laplace operators and Sobolev
calculus in spaces of differential forms, harmonic forms, Hodge
decomposition), has been a very active topic of research in recent
years. It turns out that many important examples of
infinite-dimensional nonflat spaces (loop spaces, product
manifolds, configuration spaces) are naturally equipped with
probability measures (Brownian bridge, Poisson measures, Gibbs
measures). Properties of these measures depend in a nontrivial way
on the differential geometry of the underlying spaces themselves,
and play therefore a significant role in their study. Moreover, in
many cases the absence of a proper smooth manifold structure makes
it more natural to work with $L^2$-objects (such as functions,
sections, etc.)\ on these infinite-dimensional spaces, rather than
to define analogs of the smooth ones.

Thus, the concept of an $L^2$-deRham complex has an important
meaning in this framework. The study of $L^2$-cohomologies for
finite-dimensional manifolds, initiated in \cite{Ati}, has been a
subject of many works, see e.g.\   \cite{BrLes,Dod,ELR} and  the
review papers \cite{Mat,Pan}. In the infinite-dimensional case,
loop spaces have been most studied \cite{EL,JL, LRo,  Le}, the
papers \cite{EL, Le} containing also a review of the subject. The
deRham complex on infinite product manifolds with Gibbs measures
(which appear in connection with problems of classical statistical
mechanics) was constructed in \cite{ADK1, ADK2} (see also
\cite{LeBe} for the case of the infinite-dimensional torus). We
should also mention the papers \cite{AK, Ar1, Ar2, ArM, Shi},
where the
case of a flat (Hilbert) state space has been considered (the $L^2$%
-cohomological structure turns out to be nontrivial even in this
case due to the existence of interesting measures on such a
space).

In \cite{ADL1, ADL2}, the authors started the study of
differential forms over the infinite-dimensional space $\Gamma _X$
and the corresponding Laplacians (of Bochner and deRham type)
acting in the $L^2$-spaces
with respect to a Poisson measure. In \cite{ADL3}, the associated $L^2$%
-cohomologies have been investigated.

Another approach to the construction of differential forms and
related objects over Poisson spaces, based on the ``transfer
principle'' from Wiener spaces, was proposed in \cite{Pr2}, see
also \cite{PPr} and \cite{Pr}.

It should be stressed that the choice of an underlying measure
plays a crucial role in all these studies. The results  of
\cite{ADL1,ADL2,ADL3} have only covered the case of  Poisson
measures, which are related to mathematical models of ``free''
systems, i.e., systems without interaction. The choice of more
complicated measures, such as Gibbs type perturbations of Poisson
measures,  is particularly motivated by the study of interacting
systems  of classical statistical mechanics. Properties of the
corresponding Laplace operators may then strongly depend on the
choice of an appropriate measure.

In order to develop a reasonable theory covering also this case,
we need to restrict ourselves to a class of measures on $\Gamma_X$
that possess a certain regularity. So, we consider those measures
$\mu$ which satisfy the
following condition:  for any measurable function $F:\Gamma _X\times X\to {\Bbb R}$, $%
F\ge 0$,
\begin{equation}
\int_{\Gamma _X}\mu (d\gamma )\,\sum_{x\in \gamma }F(\gamma
,x)=\int_{\Gamma _X}\mu (d\gamma )\int_X\sigma (\gamma
,dx)\,F(\gamma \cup\{x\},x), \label{mecke1}
\end{equation}
where $\sigma (\gamma ,\cdot)$ is a Borel measure on $X$ which is
absolutely continuous with respect to the volume measure $m$ on
$X$ for $\mu $-a.e.\ $\gamma \in \Gamma _X$. In particular, the
Poisson measure with intensity $\rho(x)\, m(dx)$ satisfies
\eqref{mecke1} with $\sigma(\gamma,dx)=\rho(x)\,m(dx)$, and in
this case \eqref{mecke1} becomes the classical Mecke identity
\cite{Mec67}, see also \cite{Ka75,KMM}. Furthermore, as shown by
Georgii \cite{Ge79} and  Nguyen and Zessin \cite{NZ},
\eqref{mecke1} holds for all Gibbs measures. The class of all
probability measures on $\Gamma_X$ satisfying \eqref{mecke1} was
singled out in \cite{MWM} (see also \cite{WE}), where
\eqref{mecke1} was called condition $\Sigma_m'$. A relation
between this condition  and an integration by parts formula for a
measure $\mu$ was studied in \cite{Lipsche}.

An iterated application of \eqref{mecke1} to a function
$F:\Gamma_X\times X^k\to{\Bbb R}$, $k\in\N$, gives rise to a
family of random measures $\sigma^{(k)}(\gamma)$ on $X^k$.

The structure of the present paper is as follows. In Section~2 we
recall the definition of a differential form over $\Gamma _X$,
first given in \cite{ADL1,ADL2}, and introduce the spaces $L_\mu
^2\Omega ^n$ of forms that are square integrable with respect to
$\mu $. We construct a unitary isomorhism
\begin{equation}\label{awq9876} I^n:L_\mu ^2\Omega ^n\rightarrow
\bigoplus_{k=1}^n L_\mu ^2\big(\Gamma _X\to \bigcup\limits_{\gamma
\in \Gamma_X} L_{\sigma ^{(k)}(\gamma )}^2\Psi _{\mathrm
sym}^n(X^k)\big), \end{equation} where $L_\mu ^2\big(\Gamma _X\to
\bigcup\limits_{\gamma \in \Gamma_X} L_{\sigma
^{(k)}(\gamma )}^2\Psi _{\mathrm sym}^n(X^k)\big)$ is the space of $\mu $%
-square-integrable mappings \begin{equation}\label{332633} \Gamma
_X\ni \gamma \mapsto W(\gamma )\in L_{\sigma ^{(k)}(\gamma
)}^2\Psi _{\mathrm sym}^n(X^k), \end{equation} and $L_{\sigma
^{(k)}(\gamma )}^2\Psi _{ \mathrm sym}^n(X^k)$ is a space of
$n$-forms over $X^k$ that are square-integrable with respect to
$\sigma ^{(k)}(\gamma )$ and satisfy some additional conditions.
In the case where $\mu $ is a  Poisson measure $\pi$, the
isomorphism $I^n$ was constructed in \cite{ADL3}.

In Section 3, we define  Bochner type operators in $L_\mu ^2\Omega
^n$. First, we introduce the bilinear form $$ {\cal E}_{\mu,n}
^{\mathrm B}(W^{(1)},W^{(2)}):=\int_{\Gamma
_X}\langle\nabla^\Gamma W^{(1)}(\gamma ),\nabla^\Gamma
W^{(2)}(\gamma )\rangle\,\mu (d\gamma ) $$ on the space of smooth
local forms, where $\nabla^\Gamma $ is the covariant derivative on
$\Gamma _X$ (introduced in \cite{ADL1, ADL2}), and prove its
closability. We call the corresponding generator ${\bf H}_{\mu,n}
^{\mathrm B}$ the Bochner Laplacian on $\Gamma _X$ associated with
$\mu $.

Further, we show that,  under the action of the isomorphism $I^n$,
the form ${\cal E}_{\mu,n} ^{\mathrm B}$ can  be expressed via
Bochner type bilinear forms ${\cal E}_{\sigma ^{(k)}(\gamma
)}^{\mathrm B}$ associated with the measures $\sigma ^{(k)}(\gamma
)$ on $X^k$, $k=1,\dots,n$, $\mu$-a.e.\ $\gamma \in \Gamma _X $.
As an application of this result, we derive
sufficient conditions for the space of all Bochner-harmonic $\mu $%
-square integrable forms on $\Gamma _X$ to consist only of the
zero form. Let us remark that we do not assume extremality of
$\mu$, so that  nonconstant $\mu $-square integrable harmonic
functions on $\Gamma _X$ may in general exist \cite{AKR2}.

In Section 4, we introduce and study the structure of the deRham
complex in the spaces $L_\mu ^2\Omega ^n$. Following \cite{ADL3},
we first define
 a Hodge--deRham differential ${\bf d}_n$ on
the space of smooth local forms. We prove the closability of the
${\bf d}_n$'s as operators from $L_\mu ^2\Omega ^n $ into $ L_\mu
^2\Omega ^{n+1}$ and consider the Hilbert complex $$
\cdots \stackrel{\bar{{\bf d}}_{n-1}}{\longrightarrow }L_\mu ^2\Omega ^n%
\stackrel{\bar{{\bf d}}_n}{\longrightarrow }L_\mu ^2\Omega ^{n+1}\stackrel{%
\bar{{\bf d}}_{n+1}}{\longrightarrow }\cdots \,, $$ where ${\bf
\bar{d}}_n$'s are the corresponding closures. Next, we define a
Hodge--deRham Laplacian ${\bf H}_{\mu,n} ^{\mathrm R}$ as the
generator of the closed form $$
{\cal E}_{\mu,n} ^{\mathrm R}(W^{(1)},W^{(2)}):=( {\bf \bar{d}}_nW^{(1)},{\bf \bar{d}}%
_nW^{(2)}) _{L_\mu ^2\Omega ^{n+1}}+( {\bf d}_{n-1}^{*}W^{(1)},{\bf d}%
_{n-1}^{*}W^{(2)}) _{L_\mu ^2\Omega ^{n-1}} $$ on $L_\mu ^2\Omega
^n$ with  domain $D({\cal E}_{\mu,n} ^{\mathrm R})=D(\bar{\bf
d}_n)\cap D({\bf d}_{n-1}^{*})$. We prove that, under certain
additional conditions on $\mu $, the domain of the operator ${\bf
H}_{\mu,n} ^{\mathrm  R}$ contains smooth local forms. This gives
us a possibility to prove, for ${\bf H}_{\mu,n}^{\mathrm B}$ and
${\bf H}_{\mu,n}^{\mathrm R}$, an analog of the Weitzeb\"{o}ck
formula.

In Section 5, we consider our main examples: Gibbs measures with
pair interaction on $\Gamma_{X}$. More exactly, we consider in
details Ruelle type measures on $\Gamma_{\R^d}$ (cf.\
\cite{Ru70}), and Gibbs measures in the low activity--high
temperature regime, as well as Gibbs measures with positive
potentials on $\Gamma_X$.
 In these cases,
we get more explicit expressions for the Bochner and deRham
Laplacians.

It is a great pleasure to thank D.~ B.~Applebaum, K.~D.~Elworthy,
 P.~Malliavin, M.~R\"{o}ckner for their interest
in this work and helpful discussions. The financial support of SFB
256, DFG Research Projects 436 RUS 113/593 and 436 UKR 113/43,
BMBF Research Project UKR-004-99, and the British-German research
project 313/ARC-pz-XIII-99 is gratefully acknowledged.

\section{Differential forms over a configuration space}
\label{ewew086}

Let $X$ be a complete,  connected, oriented, $C^\infty $
Riemannian manifold of infinite volume. Let $d$ denote the
dimension of $X$. Let $\langle \cdot ,\cdot
\rangle _x$ denote the  inner product in the tangent space $T_xX$ to $%
X $ at a point $x\in X$. The associated norm will be denoted by
$|\cdot |_x $. Let  $\nabla ^X$ stand for the gradient on $X$.

The configuration space $\Gamma _X$ over $X$ is defined as the set
of all locally finite subsets (configurations) in $X$:
\begin{equation*}
\Gamma _X:=\big\{ \,\gamma \subset X\mid |\gamma_ \Lambda |<\infty
\text{ for each compact }\Lambda \subset X\,\big\} . \nonumber
\end{equation*}
Here, $\gamma_\Lambda{:=\gamma\cap\Lambda}$ and  $|A|$ denotes the
cardinality of a set $A$.

We can identify any $\gamma \in \Gamma _X$ with the positive,
integer-valued Radon measure
\begin{equation*}
\sum_{x\in \gamma }\varepsilon _x\in {\cal M}(X),  \nonumber
\end{equation*}
where $\varepsilon _x$ is the Dirac measure with mass at $x$,
$\sum_{x\in \varnothing }\varepsilon _x:=$zero measure, and ${\cal
M}(X)$ denotes the set of all positive Radon measures on the Borel
$\sigma $-algebra ${\cal B} (X)$. The space $\Gamma _X$ is endowed
with the relative topology as a subset of the space ${\cal M}(X)$
with the vague topology, i.e., the weakest topology on $\Gamma _X$
with respect to which  all maps
\begin{equation*}
\Gamma _X\ni \gamma \mapsto \langle f,\gamma \rangle
:=\int_Xf(x)\,\gamma (dx)\equiv \sum_{x\in \gamma }f(x)  \nonumber
\end{equation*}
are continuous. Here, $f\in C_0(X)$($:=$the set of all continuous
functions on $X$ with compact support). Let ${\cal B}(\Gamma _X)$
denote the corresponding Borel $\sigma $-algebra.

The tangent space to $\Gamma _X$ at a point $\gamma $ is defined
as the Hilbert space
\begin{equation}
T_\gamma \Gamma _X{:=}L^2(X\to TX;\gamma)\equiv\bigoplus_{x\in
\gamma }T_xX. \label{tg-sp1}
\end{equation}
The scalar product and the norm in $T_\gamma \Gamma _X$ will be denoted by $%
\langle \cdot ,\cdot \rangle _\gamma $ and $\left\| \cdot \right\|
_\gamma $, respectively. Thus, each $V(\gamma )\in T_\gamma \Gamma
_X$ has the form $V(\gamma )=(V(\gamma, x))_{x\in \gamma }$, where
$V(\gamma, x)\in T_xX$, and
\begin{equation*}
\| V(\gamma )\| _\gamma ^2=\sum_{x\in \gamma }|V(\gamma ,x)|_x^2.
\nonumber
\end{equation*}

We now recall how to define derivatives of a function
$F:\Gamma_X\to\R$. Let $\gamma \in \Gamma _X$ and $x\in \gamma $.
By ${\cal O}_{\gamma ,x}$ we  denote an arbitrary open
neighborhood of $x$ in $X$ such that ${\cal O}_{\gamma
,x}\cap(\gamma \setminus \{x\})=\varnothing$. We define the
function $${\cal O}_{\gamma,x}\ni y\mapsto
F_x(\gamma,y){:=}F(\gamma-\eps_x+\eps_y)\in\R.$$ We say that $F$
is differentiable at $\gamma\in\Gamma_X$ if, for each
$x\in\gamma$, the function $F_x(\gamma,\cdot)$ is differentiable
at $x$ and $$\nabla^\Gamma F(\gamma){:=}(\nabla^X_x F
(\gamma))_{x\in\gamma}\in T_\gamma\Gamma_X,\qquad \nabla^X_x
F(\gamma){:=}\nabla^X F_x(\gamma,x). $$ Analogously, the higher
order derivatives of $F$ are defined,
$(\nabla^\Gamma)^{(k)}F(\gamma)\in(T_\gamma\Gamma_X)^{\otimes k}$,
$k\in\N$.

Let ${\cal O}_c(X)$ denote the set of all open relatively compact
sets in $X$. A function $F:\Gamma_X\to\R$ is called local if there
exists  $\Lambda\in{\cal O}_c(X)$ such that
$F(\gamma)=F(\gamma_\Lambda)$ for each $\gamma\in\Gamma_X$.

Any function of the form
\begin{equation}\label{errdtftfz}
F(\gamma)=g_F(\langle\varphi_1,\gamma\rangle,\dots,\langle\varphi_N,\gamma\rangle),
\end{equation}
where $g_F\in C^\infty_{\mathrm b}(\R^N)$ and
$\varphi_1,\dots,\varphi_N\in {\cal D}:=C_0^\infty(X)$(:$=$the set
of all infinitely differentiable functions on $X$ with compact
support), is local, bounded, infinitely differentiable, and the
derivatives of $F$ are polynomially bounded:
\begin{equation}\label{553434}\forall k\in\N\ \exists \varphi\in C_0(X),\
\varphi\ge0:\
\|(\nabla^\Gamma)^{(k)}F(\gamma)\|^2_{(T_\gamma\Gamma_X)^{\otimes
k}} \le \langle \varphi,\gamma\rangle^k\qquad\text{for all
}\gamma\in\Gamma_X.\end{equation} The set of all functions of the
form \eqref{errdtftfz} will be denoted by ${\cal
FC}^\infty_{\mathrm b}({\cal D},\Gamma_X)$.

Vector fields and first order differential forms on $\Gamma _X$
will be identified with sections of the bundle $T\Gamma _X$.
Higher order
differential forms will be identified with sections of the tensor bundles $%
\wedge ^n(T\Gamma _X)$ with fibers
\begin{equation}
\wedge ^n(T_\gamma \Gamma _X)%
=\wedge ^n\left( \bigoplus_{x\in \gamma }T_xX%
\right),   \label{tang-n}
\end{equation}
where $\wedge ^n({\cal H})$ (or ${\cal H}^{\wedge n}$) stands for
the $n$th antisymmetric tensor power of a Hilbert space ${\cal
H}$.
Thus, under a differential form $W$ of order $n$, $n\in {\Bbb N}$, over $%
\Gamma _X,$ we will understand a mapping
\begin{equation}
\Gamma _X\ni \gamma \mapsto W(\gamma )\in \wedge ^n(T_\gamma
\Gamma _X). \label{lkghf}
\end{equation}

We will now recall how to introduce a  covariant derivative of a
differential form \eqref{lkghf}.

Let $\gamma\in\Gamma_X$ and $x\in\gamma$. We  define  the mapping
\begin{equation*}\notag
{\cal O}_{\gamma ,x}\ni y\mapsto W_x(\gamma ,y)%
\mbox{$:=$}%
W(\gamma _y)\in \wedge ^n(T_{\gamma _y}\Gamma _X),\qquad
\gamma_y:=\gamma-\eps_x+\eps_y. \label{sec1}
\end{equation*}
This  is a section of the Hilbert bundle
\begin{equation}
\wedge ^n(T_{\gamma _y}\Gamma _X)\mapsto y\in {\cal O}_{\gamma
,x}. \label{bund1}
\end{equation}
The Levi--Civita connection on $TX$ generates in a natural way a
connection on this bundle. We denote by $\nabla _{\gamma ,x}^X$
the corresponding covariant derivative and use the notation
\begin{equation*}
\nabla _x^XW(\gamma )%
\mbox{$:=$}%
\nabla _{\gamma ,x}^X\,W_x(\gamma ,x)\in T_xX\otimes \left( \wedge
^n(T_\gamma \Gamma _X)\right)  \nonumber
\end{equation*}
if the section $W_x(\gamma ,\cdot )$ is differentiable at $x$.

We say that the form $W$ is differentiable at a point $\gamma $ if
for each $x\in \gamma $ the section $W_x(\gamma ,\cdot )$ is
differentiable at $x$, and
\[
\nabla ^\Gamma W(\gamma )%
\mbox{$:=$}%
(\nabla _x^XW(\gamma ))_{x\in \gamma }\in T_\gamma \Gamma
_X\otimes \left( \wedge ^n(T_\gamma \Gamma _X)\right) .
\]
The mapping {\
\begin{equation*}\notag
\Gamma _X\ni \gamma \mapsto \nabla ^\Gamma W(\gamma )%
\in T_\gamma \Gamma _X\otimes \left( \wedge ^n(T_\gamma \Gamma
_X)\right)
\end{equation*}
will be called the covariant gradient of the form $W$.

Analogously, one can introduce higher order derivatives of a
differential form $W$. Precisely, the $k$th derivative
$(\nabla^\Gamma)^{(k)}W(\gamma)$ belongs to
$(T_\gamma\Gamma_X)^{\otimes k}\otimes
(\wedge^n(T_\gamma\Gamma_X))$.

Let us note that, for any $\eta \subset \gamma $, the space
$\wedge
^n(T_\eta \Gamma _X)$ can be identified in a natural way with a subspace of $%
\wedge ^n(T_\gamma \Gamma _X)$. In this sense, we will use the
expression  $W (\gamma )=W (\eta )$ without additional
explanations.

A form $W:\Gamma_X\to\wedge^n(T\Gamma_X)$ is called local  if
there exists  $\Lambda=\Lambda(W)\in {\cal O}_c(X)$  such that
$W(\gamma)=W(\gamma_\Lambda)$ for each $\gamma\in\Gamma_X$.

Let ${\cal F}\Omega^n$ denote the set of all local, infinitely
differentiable forms $W:\Gamma_X\to\wedge ^n(T\Gamma_X)$ such that
there exist  $\varphi\in C_0(X)$, $\varphi\ge0$, and $l\in\N$
(depending on $W$) satisfying:
\begin{equation}\label{weghtse}\|W(\gamma)\|^2_{\wedge^n(T_\gamma\Gamma_X)}
\le\langle\varphi,\gamma\rangle^l\qquad\text{for all
}\gamma\in\Gamma_X.\end{equation}
 Below, we will give an explicit
construction of a class of forms belonging to ${\cal F}\Omega^n$.

Let $\mu$ be a probability measure on $(\Gamma_X,{\cal
B}(\Gamma_X))$ which has all moments finite, i.e.,
\begin{equation}\label{edrdt}\forall k\in\N,\ \forall \varphi\in C_0(X),\
\varphi\ge0:\qquad
\int_{\Gamma_X}\la\varphi,\gamma\ra^k\,\mu(d\gamma)<\infty.\end{equation}

Our next goal is to give a description of the space of $n$-forms
that are square-integrable with respect to the measure $\mu$.

Let $\widetilde{{\cal F}\Omega ^n}^\mu$ denote the $\mu$-classes
determined by ${\cal F}\Omega ^n$.  We define on $\widetilde{{\cal
F}\Omega ^n}^\mu$ the $L^2$-scalar product with respect to the
measure $\mu$:
\begin{equation}
(W_1,W_2)_{L_{\mu }^2\Omega ^n}%
\mbox{$:=$}%
\int_{\Gamma _X}\langle W_1(\gamma ),W_2(\gamma )\rangle _{\wedge
^n(T_\gamma \Gamma _X)}\,\mu (d\gamma ).  \label{4.1}
\end{equation}
The integral on the right hand side of \eqref{4.1} is finite
because of \eqref{weghtse} and \eqref{edrdt}. Now, we define the
Hilbert space
$
L_{\mu }^2\Omega ^n%
=L^2(\Gamma _X\to \wedge ^n(T\Gamma _X);\mu )  
$ as the completion of $\widetilde{{\cal F}\Omega ^n}^\mu$ with
respect to the norm generated by the scalar product (\ref{4.1}).
In what follows, we will not distinguish in notations between
${\cal F}\Omega ^n$ and $\widetilde{{\cal F}\Omega ^n}^\mu$, since
it will be clear from the context which of these sets we mean.

Let $m$ denote the volume measure on $X$. From now on, we suppose
that, for any measurable function $F:\Gamma_X\times X\to\R$,
$F\ge0$,
\begin{equation}\label{qwgfhjzjg}
\int_{\Gamma_X}\mu(d\gamma)\int_X \gamma(dx)\, F(\gamma
,x)=\int_{\Gamma_X}\mu(d\gamma)\int_X\sigma (\gamma,dx)\,
F(\gamma+\eps_x,x) ,\end{equation} where $\sigma(\gamma,\cdot)\ll
m$ for $\mu$-a.e.\ $\gamma\in\Gamma_X$. We shall use the notation
\begin{equation*}\label{esesa}\rho(\gamma,x){:=}
\frac{d\sigma(\gamma,\cdot)}{dm}(x).\end{equation*} In the theory
of point processes,  this property of the measure $\mu$ is called
$\Sigma_m'$, see \cite{MWM}. All Gibbs measures, in particular,
all Poisson measures satisfy this property, see \cite{Ge79,Mec67,
NZ}. We consider this case in Section~\ref{dseersre}.

We will need the following consequence of the property
$\Sigma_m'$. Let ${:}\,\gamma ^{\otimes k}\,{:}$ be the measure on
$X^k$ given by
\begin{equation*}
{:}\,\gamma ^{\otimes k}\,{:}(dx_1,\dots ,dx_k):=\sum_{\left\{
y_1,\dots
,y_k\right\} \subset \gamma }\varepsilon _{y_1}\widehat{\otimes }\dotsm%
\widehat{\otimes }\varepsilon _{y_k}(dx_1,\dots ,dx_k),  \nonumber
\end{equation*}
where
\[
\varepsilon _{y_1}\widehat{\otimes }\dotsm\widehat{\otimes
}\varepsilon _{y_k}(dx_1,\dots ,dx_k):=\frac 1{k!}\sum_{\sigma \in
S_k}\varepsilon _{y_{\sigma (1)}}\otimes \dots \otimes \varepsilon
_{y_{\sigma (k)}}(dx_1,\dots ,dx_k),
\] $S_k$ denoting the group of all permutations of
$\{1,\dots,k\}$.

For $\mu$-a.e.\ $\gamma\in\Gamma_X$, we denote by
$\sigma^{(k)}(\gamma,\cdot)$ the measure on $X^k$ given by
$$\sigma^{(k)}(\gamma,dx_1,\dots,dx_k){:=}\sigma(\gamma,dx_1)\sigma(\gamma+\eps_{x_1},
dx_2)\dotsm \sigma(\gamma+\eps_{x_1}+\dots+\eps_{x_{k-1}},dx_k),$$
and let $\mu^{(k)}$ be the measure on $\Gamma_X\times X^k$ defined
by
$$\mu^{(k)}(d\gamma,dx_1,\dots,dx_k){:=}\mu(d\gamma)\sigma^{(k)}(\gamma,dx_1,\dots,dx_k).$$

\begin{lemma}\label{waedfdrtdtgfz} For any measurable $F:\Gamma_X\times X^k\to\R$\rom,
$F\ge 0$\rom, $k\in\N$\rom,
\begin{multline}\label{awrdrt}
k!\int_{\Gamma_X}\mu(d\gamma)\int_{X^k}{:}\,\gamma ^{\otimes
k}\,{:}(dx_1,\dots ,dx_k) \, F(\gamma,x_1,\dots,x_k)
\\ =\int_{\Gamma_X\times X^k}\mu^{(k)}(d\gamma,dx_1,\dots,dx_k)\, F(\gamma+\eps_{x_1}+\dots+\eps_{x_k},
x_1,\dots,x_k).
\end{multline}\end{lemma}

\noindent {\it Proof}. We prove this by induction. For $k=1$,
\eqref{awrdrt} is just \eqref{qwgfhjzjg}. Let us suppose that
\eqref{awrdrt} holds up to $k-1$. As easily seen, $$k
\,{:}\,\gamma ^{\otimes k}\,{:}(dx_1,\dots
,dx_k)=\gamma(dx_k){:}\,(\gamma-\eps_{x_k} )^{\otimes
k-1}\,{:}(dx_1,\dots ,dx_{k-1}).$$ Then, by
 the induction hypothesis we have
\begin{gather*}k!\int_{\Gamma_X}\mu(d\gamma)\int_{X^k}{:}\,\gamma
^{\otimes k}\,{:}(dx_1,\dots ,dx_k) \, F(\gamma,x_1,\dots,x_k)
\\ =\int_{\Gamma_X}\mu(d\gamma)\int_X\gamma(dx_k)\,(k-1)! \int_{X^{k-1}}
{:}\,(\gamma-\eps_{x_k}) ^{\otimes (k-1)}\,{:}(dx_1,\dots
,dx_{k-1})F(\gamma,x_1,\dots,x_k)\\
=\int_{\Gamma_X}\mu(d\gamma)\int_X
\sigma(\gamma,dx_k)\int_{X^{k-1}}{:}\,\gamma ^{\otimes
(k-1)}\,{:}(dx_1,\dots
,dx_{k-1})F(\gamma+\eps_{x_k},x_1,\dots,x_k)\\
=\int_{\Gamma_X}\mu(d\gamma)\int_{X^{k-1}}{:}\,\gamma ^{\otimes
k-1}\,{:}(dx_1,\dots ,dx_{k-1})\int_X
\sigma(\gamma,dx_k)F(\gamma+\eps_{x_k},x_1,\dots,x_k)\\
=\int_{\Gamma_X}\mu(d\gamma)\int_X\sigma(\gamma,dx_1)\int_X\sigma(\gamma+\eps_{x_1},
dx_2)\dotsm\\ \dotsm \int_X \sigma(\gamma+\eps_{x_1}+\dots+ \eps_{
x_{k-1}},dx_k)\,F(\gamma+\eps_{x_1}+\dots+\eps_{x_k},x_1,\dots,x_k).\quad
\blacksquare
\end{gather*}

We will now give an isomorphic description of the space
$L^2_\mu\Omega^n$.  We first need some preparations. Let
$$\widetilde X^k:=\left\{(x_1,\dots,x_k)\in X^k:\ x_i\ne x_j\text{
if }i\ne j\right\}.$$ Notice that the set $X^k\setminus \widetilde
X^k$ is of zero $m^{\otimes k}$ measure.
 We have, for each $(x_1,\dots,x_k)\in\widetilde X^k$,
\begin{equation}
\wedge ^n(T_{(x_1,\dots ,x_k)}X^k)=\wedge^n\left(\bigoplus_{i=1}^k
T_{x_i}X \right) =\bigoplus_{\begin{gathered}{\scriptstyle{
0\le l_1,\dots,l_k\le d}} \\ \scriptstyle l_1+\dots+l_k=n \end{gathered}%
}(T_{x_1}X)^{\wedge l_1}\wedge \dots \wedge (T_{x_k}X)^{\wedge
l_k}. \label{n-forms}
\end{equation}
For a form $\omega : X^k\to\wedge^n(T X^k)$ and
$(x_1,\dots,x_k)\in\widetilde X^k$, we denote by $\omega
(x_1,\dots ,x_k)_{l_1,\dots ,l_k}$ the corresponding component of
$\omega (x_1,\dots ,x_k)$ in the decomposition (\ref{n-forms}).

We introduce a set $\Psi _{\mathrm sym}^n(X^k)$  of smooth forms
$\omega :X^k\to\wedge^n(TX^k)$  which have compact support and
satisfy on $\widetilde X^k$ the following assumptions:

\begin{enumerate}
\item[(i)]  $\omega (x_1,\dots ,x_k)_{l_1,\dots ,l_k}=0$ if $l_j=0$ for some
$j\in \{1,\dots ,k\}$;

\item[(ii)]  $\omega $ is invariant under the action of the group $S_k$:
\begin{equation}
\omega (x_1,\dots ,x_k)=\omega (x_{\sigma (1)},\dots ,x_{\sigma
(k)})\qquad \text{for each }\sigma \in S_k.  \label{symmetric}
\end{equation}
(we identify the spaces $T_{(x_1,\dots ,x_k)}X^k=\bigoplus_{i=1}^k
T_{x_i}X $ and $T_{(x_{\sigma (1)},\dots ,x_{\sigma
(k)})}X^k=\bigoplus_{i=1}^k T_{x_{\sigma(i)}}$ through the natural
isomorphism).
\end{enumerate}

Using \eqref{edrdt} and Lemma~\ref{waedfdrtdtgfz}, we easily
conclude that any mapping of the form
\begin{equation}\label{serfzdtcf} \Gamma_X\times X^k\ni(\gamma,
x_1,\dots,x_k)\mapsto
F(\gamma)\omega(x_1,\dots,x_k)\in\wedge^n(T_{(x_1,\dots,x_k)}X^k),\end{equation}
where $F\in{\cal F}C_{\mathrm b}^\infty({\cal D},\Gamma_X)$ and
$\omega\in\Psi_{\mathrm sym}^n(X^k)$ belongs to the space
$L^2(\Gamma_X\times X^k\to \wedge^n(TX^k);\mu^{(k)})$. Let
$L^2_\Psi(\Gamma_X\times X^k\to \wedge^n(TX^k) ;\mu^{(k)})$ denote
the closed linear span of all mappings of the form
\eqref{serfzdtcf} in $L^2(\Gamma_X\times X^k\to
\wedge^n(TX^k);\mu^{(k)})$. It is not hard to show that the latter
is just the space of all $\mu^{(k)}$-square integrable mappings of
the form $$\Gamma_X\times \widetilde
X^k\ni(\gamma,x_1,\dots,x_k)\mapsto{\cal W}(\gamma,x_1,\dots,x_k)
\in {\Bbb T}_{\{x_1,\dots ,x_k\}}^{(n)}X^k$$ such that, for
$\mu^{(k)}$-a.e.\ $(\gamma,x_1,\dots,x_k)\in\Gamma_X\times
\widetilde X^k$, $${\cal W}(\gamma,x_1,\dots,x_k)={\cal
W}(\gamma,x_{\sigma(1)},\dots,x_{\sigma(k)}),\qquad \sigma\in
S_k.$$ Here, \begin{equation}\label{eww4a}{\Bbb T}_{\{x_1,\dots
,x_k\}}^{(n)}X^k:
=\bigoplus_{\begin{gathered}{%
\scriptstyle{ 1\le l_1,\dots,l_k\le d}} \\ \scriptstyle
l_1+\dots+l_k=n
\end{gathered}}(T_{x_1}X)^{\wedge l_1}\wedge \dots \wedge (T_{x_k}X)^{\wedge
l_k}.\end{equation} (Notice that the space ${\Bbb T}_{\{x_1,\dots
,x_k\}}^{(n)}X^k$ is indeed independent of the order of the points
$x_1,\dots,x_k$.)

\setcounter{remark}{1}

\begin{remark}\rom{ Evidently, $$L^2_\Psi(\Gamma_X\times X^k\mapsto \wedge^n(TX^k)
;\mu^{(k)})=L^2_\mu
\big(\Gamma_X\to\bigcup_{\gamma\in\Gamma_X}L^2_{\sigma^{(m)}(\gamma)}\Psi_{\mathrm
sym }^n(X^m)\big),$$ where the latter space was defined in the
Introduction (see formulas \eqref{awq9876},
\eqref{332633}).}\end{remark}

By virtue of \eqref{tang-n} and \eqref{eww4a}, we have
\begin{equation} \wedge ^n(T_\gamma \Gamma
_X)=\bigoplus_{k=1}^n\bigoplus_{\{x_1,\dots ,x_k\}\subset \gamma
}{\Bbb T}_{\{x_1,\dots ,x_k\}}^{(n)}X^k. \label{tang-n1}
\end{equation}
For $W\in \Gamma_X\to \wedge^n(T\Gamma_X)$, we denote by
$W_k(\gamma )\in \bigoplus_{\{x_1,\dots ,x_k\}\subset \gamma
}{\Bbb T}_{\{x_1,\dots ,x_k\}}^{( n)}X^k$ the corresponding
component of $W(\gamma )\in \wedge ^n(T_\gamma \Gamma _X)$
 in the decomposition %
\eqref{tang-n1}. Thus, for $\{x_1,\dots ,x_k\}\subset \gamma $,
$W_k(\gamma,x_1,\dots ,x_k)$ is equal to the projection of
$W(\gamma )$ onto the subspace ${\Bbb T}_{\{x_1,\dots
,x_k\}}^{(n)}X^k$.

\setcounter{proposition}{2 }
\begin{proposition}
\label{sq-int}The space $L_{\mu }^2\Omega ^n$ is unitarily
isomorphic to the space
\begin{equation}
\bigoplus_{k=1}^n L^2_\Psi(\Gamma_X\times X^k\to\wedge
^n(TX^k);\mu^{(k)}), \label{tensor-n}
\end{equation}
where the corresponding isomorphism $I^n$ is defined by the
formula
\begin{equation}
I_k^nW(\gamma ,x_1,\dots,x_k)%
\mbox{$:=$}%
(k!)^{-1/2}\,W_k(\gamma +\eps_{x_1}+\dots+\eps_{x_k},x_1,\dots,x_k
),\qquad k=1,\dots ,n. \label{asfdgfsdd}
\end{equation}
Here{\rm ,} $I_k^nW:=\left( I^nW\right) _k$ is the $k$-th component of $%
I^nW $ in the decomposition {\rm (\ref{tensor-n}).}
\end{proposition}

\noindent{\it Proof}. A direct calculation shows that
\begin{equation}
\|W(\gamma)\|^2_{\wedge ^n(T_\gamma \Gamma
_X)}=\sum_{k=1}^n\int_{X^k}\| W_k( \gamma,x_1,\dots,x_k) \|^2
_{{\Bbb T}_{\left\{x_1,\dots,x_k\right\} }^{(n)}X^k}{:}\,\gamma
^{\otimes k}\,{:}\left( dx_1,\dots,dx_k\right). \label{tensor-ng}
\end{equation}
Therefore, by Lemma~\ref{waedfdrtdtgfz}, we have for any $W\in
{\cal F}\Omega^n$
\begin{align*}&\int_{\Gamma_X}\|W(\gamma)\|^2_{\wedge ^n(T_\gamma \Gamma
_X)}\,\mu(d\gamma)\\&\qquad=\sum_{k=1}^n\int_{\Gamma_X\times X^k
}\| W_k( \gamma+\eps_{x_1}+\dots+\eps_{x_k},x_1,\dots,x_k) \|^2
_{{\Bbb T}_{\left\{x_1,\dots,x_k\right\}
}^{(n)}X^k}\,\mu^{(k)}\left(d\gamma,
dx_1,\dots,dx_k\right).\end{align*} Hence, $I^n$ is an isometry of
the space $L^2_\mu\Omega^n$ into the space \eqref{tensor-n}. Next,
the image of each mapping \eqref{serfzdtcf} under $(I^{(n)})^{-1}$
is given by \begin{equation}\label{kikimora}
W_l(\gamma,x_1,\dots,x_l):=\begin{cases}0,& l\ne k,\\ (k!)^{1/2}
F(\gamma-  \eps_{x_1}-\dots-\eps_{x_k}
)\om(x_1,\dots,x_k),&l=k,\end{cases}\end{equation} and evidently
belongs to ${\cal F}\Omega^n$. Therefore, $I^n$ is ``onto.''\quad
$\blacksquare$\vspace{2mm}

In what follows, we will denote by ${\cal D}\Omega ^n$ the linear
span of the forms defined by (\ref{kikimora}) with $k=1,\dots,n$.
As we already noticed in the proof of Proposition \ref{sq-int},
${\cal D}\Omega ^n$ is a subset of
${\cal F} \Omega ^n$ and is dense in $L_{\mu }^2\Omega ^n$%
.

\section{Laplace operators on differential forms over configuration spaces%
}\label{dodf}

In this section, we introduce differential operators associated
with the measure $\mu$ on $\Gamma _X$ which act in the space of
square-integrable forms. These operators generalize the notions of
Bochner and deRham Laplacians on finite-dimensional manifolds. But
first, we consider the  Dirichlet operator in the space
$L^2(\Gamma_X;\mu)$.

\subsection{Dirichlet operator on functions}
\label{subsec3.1}

For each $\gamma \in \Gamma _X$, consider the triple
\begin{equation*}
T_{\gamma ,\,\infty }\Gamma _X\supset T_\gamma \Gamma _X\supset
T_{\gamma ,0}\Gamma _X.  \nonumber
\end{equation*}
Here, $T_{\gamma ,0}\Gamma _X$ consists of all finite sequences from $%
T_\gamma \Gamma _X$, and $T_{\gamma ,\,\infty }\Gamma _X%
\mbox{$:=$}%
\left( T_{\gamma ,0}\Gamma _X\right) ^{\prime }$ is the dual
space, which
consists of all sequences $V(\gamma )=(V(\gamma ,x))_{x\in \gamma }$, where $%
V(\gamma ,x)\in T_xX$. The pairing between any $V(\gamma )\in
T_{\gamma ,\,\infty }\Gamma _X$ and $v(\gamma )\in T_{\gamma
,0}\Gamma _X$ with respect to the zero space $T_\gamma \Gamma _X$
is given by
\begin{equation*}
\langle V(\gamma ),v(\gamma )\rangle _\gamma =\sum_{x\in \gamma
}\langle V(\gamma ,x),v(\gamma ,x)\rangle _x  \nonumber
\end{equation*}
(the series is, in fact, finite). From now on, under a vector field over $%
\Gamma _X$ we will understand mappings of the form $\Gamma _X\ni
\gamma \mapsto V(\gamma )\in T_{\gamma ,\infty }\Gamma _X$.

We will suppose that, for $\mu\otimes m$-a.e.\
$(\gamma,x)\in\Gamma_X\times X$, $\rho(\gamma,x)>0$ and for
$\mu$-a.e.\ $\gamma\in\Gamma_X$, the function $\rho(\gamma,\cdot)$
is weakly differentiable on $X$. We set $$\beta_\sigma(\gamma,x):=
\frac{\nabla^X_x \rho(\gamma,x)}{\rho(\gamma,x)},\qquad
\text{$\mu\otimes m$-a.e.\ $(\gamma,x)\in\Gamma_X\times X$}$$
($\beta_\sigma(\gamma,\cdot)$ is called the logarithmic derivative
of the measure $\sigma(\gamma,\cdot)$).

The logarithmic derivative of the measure $\mu$ is set to be  the
$\mu$-a.e.\ defined vector field on $\Gamma_X$ given by
\begin{gather*}\gamma\mapsto
B_\mu(\gamma)=(B_\mu(\gamma,x))_{x\in\gamma}\in
T_{\gamma,\infty}\Gamma_X,\\
B_\mu(\gamma,x):=\beta_\sigma(\gamma-\eps_x,x).\end{gather*}

We define a bilinear  form ${\cal E}_\mu$ on the space
$L^2(\Gamma_X;\mu)$ by setting
\begin{equation}\label{tfzgfzdxtfgfzgff}{\cal
E}_\mu(F^{(1)},F^{(2)}):=\int_{\Gamma_X}\langle\nabla^\Gamma
F^{(1)}(\gamma),\nabla^\Gamma
F^{(2)}(\gamma)\rangle_\gamma\,\mu(d\gamma),\end{equation} where
 $F^{(1)},F^{(2)}\in D({\cal E}_\mu){:=}{\cal F}C_{\mathrm b}^\infty ({\cal
D},\Gamma_X)$. By \eqref{553434}, \eqref{edrdt}, and
\cite[Theorem~2.4]{MR}, 
$({\cal E}_\mu,\FC)$ is a
pre-Dirichlet form.

\begin{theorem}\label{eshj} Suppose that\rom, for any  $\Lambda\in{\cal O}_c(X)$\rom,
\begin{equation}\label{568734}
\int_{\Gamma_X}\left(\sum_{x\in\gamma_\Lambda}|B_\mu(\gamma,x)|_x\right)^2\,\mu(d\gamma)<\infty.\end{equation}
Then\rom, for any $F^{(1)},F^{2)}\in {\cal F}C_{\mathrm b}^\infty
({\cal D},\Gamma_X)$\rom, we have
\begin{equation}\label{21535980}{\cal
E}_\mu(F^{(1)},F^{(2)})=\int_{\Gamma_X}({\bf H}_\mu
F^{(1)})(\gamma)F^{(2)}(\gamma)\,\mu(d\gamma),\end{equation} where
${\bf H}_\mu$ is the operator in the space $L^2(\Gamma_X;\mu)$
with domain ${\cal F}C_{\mathrm b}^\infty ({\cal D},\Gamma_X)$
given by
\begin{equation}\label{esgfsryfzsdrdfz}({\bf H}_\mu
F)(\gamma):=-\Delta^\Gamma F(\gamma)-\langle \nabla^\Gamma
F(\gamma),B_\mu(\gamma)\rangle_\gamma,\qquad F\in {\cal
F}C_{\mathrm b}^\infty ({\cal D},\Gamma_X).\end{equation} Here,
\begin{equation}\label{werweadrt}\Delta^\Gamma
F(\gamma):=\sum_{x\in\gamma}\Delta_x^X F(\gamma),\qquad \Delta_x^X
F(\gamma){:=} \Delta^X F_x(\gamma,x),\end{equation} where
$\Delta^X$ denotes the Laplacian on $X$ corresponding to the
volume measure $m.$\end{theorem}

\setcounter{corollary}{1}
\begin{corollary}
$({\cal E}_\mu,\FC)$ is closable on $L^2(\Gamma_X;\mu)$. Its
closure\rom, denoted by  $({\cal E}_\mu,D({\cal E}_\mu))$\rom,  is
associated with a positive definite self-adjoint operator, the
Friedrichs extension of ${\bf H}_\mu$\rom, which we also denote by
${\bf H}_\mu$\rom.\end{corollary}

\setcounter{remark}{2}
\begin{remark}\rom{In case of a Ruelle measure, a theorem on the $L^2$-generator of the biliear
form \eqref{tfzgfzdxtfgfzgff} was proved in \cite{AKR2}. A theorem
on the closability of the form \eqref{tfzgfzdxtfgfzgff} in the
case of a Gibbs measure on a manifold $X$ was proved in
\cite{Luis} and in the general case of a $\Sigma_m'$-measure in
\cite{MR}, see also \cite{Pedestrian}.}
\end{remark}

\noindent {\it Proof of Theorem\/} \ref{eshj}. First, we note
that, for each $F\in\FC$ and each $\gamma\in\Gamma_X$, the
function $f(x){:=}F(\gamma+\eps_x)-F(\gamma)$ belongs to $\cal D$
and $\nabla^X f(x)=\nabla^X_xF(\gamma+\eps_x)$.

Let now $F^{(1)},F^{(2)}\in\FC$ and let $\Lambda\in{\cal O}_c(X)$
be such that there exits a compact $\Lambda'\subset\Lambda$
satisfying $F^{(i)}(\gamma)=F^{(i)}(\gamma_{\Lambda'})$, $i=1,2$,
for all $\gamma\in\Gamma_X$. Then, by
\eqref{qwgfhjzjg}\begin{gather*}\int_{\Gamma_X}\la\nabla^\Gamma
F^{(1)}(\gamma),\nabla^\Gamma F^{(2)}(\gamma)\ra_\gamma\,
\mu(d\gamma)\\= \int_{\Gamma_X}\mu(d\gamma)\int_\Lambda m(dx)\,
\rho(\gamma,x) \la \nabla^X_x F^{(1)}(\gamma+\eps_x),\nabla^X_x
F^{(2)} (\gamma+\eps_x)\ra_x\\
=-\int_{\Gamma_X}\mu(d\gamma)\int_\Lambda m(dx)\,\rho(\gamma,x)
\big(\Delta_x^X F^{(1)}(\gamma+\eps_x)+\la \nabla_x^X
F^{(1)}(\gamma+\eps_x),\beta_\sigma(\gamma,x)\ra_x \big)
F^{(2)}(\gamma+\eps_x)\\ =-\int_\Gamma
\mu(d\gamma)\sum_{x\in\gamma_\Lambda} \big(\Delta_x^X
F^{(1)}(\gamma)+\la \nabla_x^X
F^{(1)}(\gamma),B_\mu(\gamma,x)\ra_x \big) F^{(2)}(\gamma)\\
=\int_\Gamma ({\bf H}_\mu
F^{(1)})(\gamma)F^{(2)}(\gamma).\end{gather*} As easily seen,
condition \eqref{568734} guarantees the inclusion ${\bf H}_\mu
F^{(1)}\in L^2(\Gamma;\mu)$.\quad $\blacksquare$\vspace{2mm}

\subsection{Bochner Laplacian on forms }

Let us consider the bilinear   form ${\cal E}_{\mu,n}^{\mathrm B}$
 defined by
\begin{equation}  \label{4.2}
{\cal E}^{\mathrm B}_{\mu,n}(W^{(1)},W^{(2)})=\int_{\Gamma_X}
\langle\nabla^\Gamma W ^{(1)}(\gamma),\nabla^\Gamma W
^{(2)}(\gamma)\rangle_{T_\gamma\Gamma_X
\otimes\wedge^n(T_\gamma\Gamma_X)}\,\mu(d\gamma),
\end{equation}
where $W^{(1)},W^{(2)}\in D({\cal E}_{\mu,n}^{\mathrm B }):={\cal
D}\Omega^n$. It follows from the definition of ${\cal D}\Omega^n$
that, for each $W\in{\cal D}\Omega^n$, there exists
$\varphi\in{\cal D}$, $\varphi\ge0$, such that
\begin{equation}\label{234539897} \|\nabla^\Gamma W(\gamma)\|^
2_{T_\gamma\Gamma_X\otimes \wedge^n(T_\gamma\Gamma_X)}\le
\la\varphi,\gamma\ra^{n+1}\qquad\text{for all
}\gamma\in\Gamma_X,\end{equation} and therefore, by \eqref{edrdt},
the function under the sign of integral in \eqref{4.2} is
integrable with respect to $\mu$.

The following lemma shows that the bilinear form $({\cal
E}_{\mu,n}^{\mathrm B},{\cal D}\Omega^n)$ is well defined on
$L^2\Omega^n$.

\setcounter{lemma}{3}

\begin{lemma} \label{uhizhgzt} We have ${\cal E}_{\mu,n}^{\mathrm
B}(W^{(1)},W^{(2)})=0$ for all $W^{(1)},W^{(2)}\in{\cal
D}\Omega^n$ such that $W^{(1)}=0$ $\mu$-a\rom.e\rom.\end{lemma}

\noindent{\it Proof}.  Let $W\in{\cal D}\Omega^n$ and $W=0$
$\mu$-a.e. For $x_0\in X$ and $R>0$, let $$B(x_0,R){:=}\{x\in
X\mid d(x_0,x)<R\},$$ where $d(\cdot,\cdot)$ denotes the
Riemannian distance on $X$. Then, \begin{align*} 0&=
\int_{\Gamma_X}\mu(d\gamma) \int_{B(x_0,R)} \gamma(dx)\,
\|W(\gamma)\|_{\wedge^n(T_\gamma\Gamma_X)}\\ &=
\int_{\Gamma_X}\mu(d\gamma)\int_{B(x_0,R)} m(dx)\, \rho(\gamma,x)
\|W(\gamma+\eps_x)\|_{\wedge^n(T_{\gamma+\eps_x}\Gamma_X)}.\end{align*}
Since $R$ was arbitrary, we therefore have
$$\|W(\gamma+\eps_x)\|_{\wedge^n(T_{\gamma+\eps_x}\Gamma_X)}=0,\qquad\text{$\mu\otimes
m$-a.e.\ $(\gamma,x)\in \Gamma_X\times X$}.$$ For a fixed
$\gamma\in\Gamma_X$, the function $X\setminus\gamma\ni x\mapsto
\|W(\gamma+\eps_x)\|_{\wedge^n(T_{\gamma+\eps_x}\Gamma_X)} $ is
continuous, and therefore for $\mu$-a.e.\ $\gamma\in\Gamma_X$,
$W(\gamma+\eps_x)=0$ on $X\setminus \gamma$. Hence, \begin{align*}
{\cal E}_{\mu,n}^{\mathrm B}(W,W)&=\int
_{\Gamma_X}\mu(d\gamma)\int_X \gamma(dx)\, \|\nabla^X_x
W(\gamma)\|_{T_x X\otimes \wedge^n(T_\gamma\Gamma_X)}^2\\ &=
\int_{\Gamma_X}\mu(d\gamma)\int_X m(dx)\,
\rho(\gamma,x)\|\nabla^X_x W(\gamma+\eps_x)\|^2_{T_x(X)\otimes
\wedge^n(T_{\gamma+\eps_x}\Gamma_X)}=0.\end{align*} From here the
lemma follows by the Schwarz inequality.\quad $\blacksquare$

\setcounter{theorem}{4}
\begin{theorem}
\label{th4.1} Suppose that
\begin{equation}\label{awrzwerjidrt}\forall\Lambda\in{\cal O}_c(X)\
 \exists\eps>0:\quad \int_{\Gamma_X}\left(\sum_{x\in\gamma_\Lambda}
|B_\mu(\gamma,x)|_x\right)^{2+\eps}\,\mu(d\gamma)<
\infty.\end{equation} Then\rom, for any $W^{(1)},W^{(2)}\in {\cal
D}\Omega ^n$, we have
\[
{\cal E}_{\mu,n }^{\mathrm B}(W^{(1)},W^{(2)})=\int_{\Gamma
_X}\langle {\bf H}_{\mu,n }^{\mathrm B}W^{(1)}(\gamma
),W^{(2)}(\gamma )\rangle _{\wedge ^n(T_\gamma \Gamma _X)}\,\mu
 (d\gamma ),
\]
where ${\bf H}_{\mu,n}^{\mathrm B}$ is the operator in the space
$L_{\mu}^2\Omega ^n$ with domain ${\cal D}\Omega ^n$ given by
\begin{equation}
{\bf H}_{\mu,n}^{\mathrm B}W(\gamma ):=-\Delta ^\Gamma W(\gamma
)-\langle \nabla ^\Gamma W(\gamma ),B_{\mu }(\gamma )\rangle
_\gamma ,\qquad W\in {\cal D}\Omega ^n.  \label{boch1}
\end{equation}
Here{\rm , }
\begin{equation}
\Delta ^\Gamma W(\gamma )%
\mbox{$:=$}%
\sum_{x\in \gamma }\Delta _x^XW(\gamma ) , \label{boch2}
\end{equation}
where $\Delta _x^X$ is the Bochner Laplacian of the bundle $\wedge
^n(T_{\gamma _y}\Gamma _X)\mapsto y\in {\cal O}_{\gamma ,x}$ with
the volume measure $m${\rm . }
\end{theorem}

\noindent{\it Proof}. We first note that, for any $W\in{\cal
D}\Omega^n$, the form ${\bf H}^{\mathrm B}_{\mu,n}W$ defined by
\eqref{boch1}, \eqref{boch2} belongs to $L^2_\mu\Omega^n$. Indeed,
as easily seen, $\Delta^\Gamma W\in{\cal F}\Omega^n$, and hence
$\Delta^\Gamma W\in L^2_\mu\Omega^n$. Next, choose any
$\Lambda\in{\cal O}_c(X)$ such that there exists a compact
$\Lambda'\subset \Lambda$ satisfying
$W(\gamma)=W(\gamma_{\Lambda'})$ for all $\gamma\in\Gamma_X$.
Then,
\begin{align}\label{tfooujipp}& \int_{\Gamma_X}\|\langle \nabla^\Gamma W(\gamma),
B_\mu(\gamma)\rangle_\gamma
\|^2_{\wedge^n(T_\gamma\Gamma_X)}\,\mu(d\gamma) \\ &\qquad =
\int_{\Gamma_X}\bigg\|\sum_{x\in\gamma_\Lambda}\langle \nabla^X_x
W(\gamma), B_\mu(\gamma,x)\rangle_x
\bigg\|^2_{\wedge^n(T_\gamma\Gamma_X)}\,\mu(d\gamma)\notag\\
&\qquad\le \int_{\Gamma_X}\left(
\sum_{x\in\gamma_\Lambda}\|\nabla^X_x
W(\gamma)\|_{T_xX\otimes\wedge^n(T_\gamma\Gamma_X)}|B_\mu(\gamma,x)|_x\right)^2\,
\mu(d\gamma).\notag
\end{align}
As easily seen, there exists $\varphi\in C_0(X)$, $\varphi\ge0$,
such that \begin{equation}\label{sewewe}\|\nabla^X_x
W(\gamma)\|_{T_xX\otimes\wedge^n(T_\gamma\Gamma_X)}^2\le
\la\varphi,\gamma\ra^n\qquad \text{for all $\gamma\in\Gamma_X$,
$x\in\gamma$.}\end{equation} Now, by using \eqref{edrdt},
\eqref{awrzwerjidrt}, \eqref{tfooujipp}, \eqref{sewewe}, and the
Schwarz inequality, we conclude that
\begin{equation}\label{tfujipp} \int_{\Gamma_X}\|\langle \nabla^\Gamma W(\gamma),
B_\mu(\gamma)\rangle_\gamma
\|^2_{\wedge^n(T_\gamma\Gamma_X)}\,\mu(d\gamma)<
\infty.\end{equation}

Next, we will need the following lemma, whose proof  follows
directly from the construction of the forms from ${\cal
D}\Omega^n$.

\setcounter{lemma}{5}
\begin{lemma}\label{dswjiu} For each fixed $W\in {\cal D}\Omega^n$
and $\gamma\in\Gamma_X$\rom, the mapping $$ X\setminus\gamma\ni
x\mapsto
\omega(x){:=}W(\gamma+\eps_x)\in\wedge^n(T_{\gamma+\eps_x}
\Gamma_X)=\wedge^n(T_\gamma\Gamma_X\oplus T_xX)$$
\rom(uniquely\rom) extends to a smooth form $$X\ni x\mapsto
\omega(x)\in\wedge^n(T_\gamma\Gamma_X\oplus T_xX),$$ and
$\nabla^X\omega=0$ on $\Lambda^c{:=}X\setminus\Lambda$, where
$\Lambda\subset X$ is compact and such that
$W(\gamma')=W(\gamma'_\Lambda)$ for all
$\gamma'\in\Gamma_X$\rom.\end{lemma}

Let  $W^{(1)},W^{(2)}\in{\cal D}\Omega^n$ and let $\Lambda\in{\cal
O}_c(X)$ be such that there exits a compact
$\Lambda'\subset\Lambda$ satisfying
$W^{(i)}(\gamma)=W^{(i)}(\gamma_{\Lambda'})$, $i=1,2$, for all
$\gamma\in\Gamma_X$. Then, by virtue of \eqref{qwgfhjzjg},
\eqref{tfujipp}, and Lemma~\ref{dswjiu} we get, analogously to the
proof of Theorem~\ref{eshj}:
\begin{align*} {\cal E}^{\mathrm
B}_{\mu,n}(W^{(1)},W^{(2)})&=\int_{\Gamma_X}
\mu(d\gamma)\int_\Lambda \gamma(dx)\,\langle \nabla^X_x
W^{(1)}(\gamma),\nabla^X_x W^{(2)}(\gamma)\rangle_{T_xX\otimes
\wedge^n(T_\gamma\Gamma_X)}\\&=\int_{\Gamma_X}\mu(d\gamma)\int_\Lambda
m(dx)\rho(\gamma,x)\, \langle\nabla^X_x
W^{(1)}(\gamma+\eps_x),\nabla^X_x
W^{(2)}(\gamma+\eps_x)\rangle_{T_xX\otimes\wedge^n(T_{\gamma+\eps_x})}\\
&=-\int_{\Gamma_X}\mu(d\gamma)\int_\Lambda
m(dx)\rho(\gamma,x)\,\big[ \langle \Delta_x^X
W^{(1)}(\gamma+\eps_x),W^{(2)}(\gamma+\eps_x)\rangle_{\wedge^n(T_{\gamma+\eps_x})}\\
&\qquad +\langle\langle\nabla^X_x
W^{(1)}(\gamma+\eps_x),\beta_\sigma(\gamma,x)\rangle_x,W^{(2)}(\gamma+\eps_x)\rangle_{
\wedge^n(T_{\gamma+\eps_x}\Gamma_X)}\big]\\&=-\int_{\Gamma_X}\mu(d\gamma)\int_\Lambda\gamma(dx)\big[
\langle \Delta^X_x
W^{(1)}(\gamma),W^{(2)}(\gamma)\rangle_{\wedge^n(T_\gamma\Gamma_X)}\\
&\qquad+\langle\langle \nabla^X_x
W^{(1)}(\gamma),B_\mu(\gamma,x)\rangle_x,W^{(2)}(\gamma)\rangle_{\wedge^n(
T_\gamma\Gamma_X)}\big]\\ &=\int_{\Gamma_X} \la{\bf H}^{\mathrm
B}_{\mu,n}W^{(1)}(\gamma),W^{(2)}(\gamma)\ra_{\wedge^n(T_\gamma\Gamma_X)}
\,\mu(d\gamma).\quad\blacksquare
\end{align*}

\setcounter{corollary}{6}
\begin{corollary} $({\cal E}^{\mathrm B}_{\mu,n},{\cal D}\Omega^n)$ is closable on $L^2_\mu\Omega^n$\rom. Its closure
$({\cal E}^{\mathrm B}_{\mu,n},D({\cal E}^{\mathrm B}_{\mu,n}))$
is associated with a positive definite\rom, self-adjoint operator,
the Friedrichs extension of  ${\bf H}^{\mathrm B}_{\mu,n}$\rom,
which we also denote by ${\bf H}^{\mathrm B}_{\mu,n}$\rom.
\end{corollary}

We define $(\widetilde{\cal E}^{\mathrm
B}_{\mu,n},D(\widetilde{\cal E}^{\mathrm B}_{\mu,n}))$ as the
image of the bilinear form $({\cal E}^{\mathrm B}_{\mu,n},D({\cal
E}^{\mathrm B}_{\mu,n}))$ under the unitary $I^n$.

\setcounter{proposition}{7}
\begin{proposition}\label{trdrtsers} Let ${\cal W}^{(1)},{\cal W}^{(2)}\in I^n \big({\cal
D}\Omega^n\big)$\rom. Then\rom, \begin{multline} \widetilde {\cal
E}_{\mu,n}^{\mathrm B}({\cal W}^{(1)},{\cal W}^{(2)})=\sum_{k=1}^n
\int_{\Gamma_X\times X^k} \mu^{(k)}(d\gamma,dx_1,\dots,dx_k)\times
\\ \times \big[ \langle \nabla^\Gamma_\gamma {\cal
W}^{(1)}(\gamma,x_1,\dots,x_k), \nabla^\Gamma_\gamma {\cal
W}^{(2)}(\gamma,x_1,\dots,x_k)\rangle_{T_\gamma\Gamma_X\otimes
{\Bbb T}^{(n)}_{\{x_1,\dots,x_k\}}X^k}\\ \text{}+ \langle
\nabla^{X^k}_{(x_1,\dots,x_k)}{\cal
W}^{(1)}(\gamma,x_1,\dots,x_k),
\nabla^{X^k}_{(x_1,\dots,x_k)}{\cal
W}^{(2)}(\gamma,x_1,\dots,x_k)\rangle_{T_{(x_1,\dots,x_k)}X^k\otimes
{\Bbb T}^{(n)}_{\{x_1,\dots,x_k\}}X^k}\big].
\label{5443}\end{multline} Here\rom, for a fixed
$(x_1,\dots,x_k)\in\widetilde X^k$\rom, $\nabla^\Gamma_\gamma$
denotes the gradient of a mapping from $\Gamma$ into ${\Bbb
T}^{(n)}_{\{x_1,\dots,x_k\}}X^k$ defined for $\mu$-a\rom.e\rom.\
$\gamma\in\Gamma_X$ similar to the gradient of a function on
$\Gamma$\rom.
\end{proposition}

\noindent {\it Proof}. Let $W^{(1)},W^{(2)}\in{\cal D}\Omega^n$
and let ${\cal W}^{(i)}{:=}I^n W^{(i)}$, $i=1,2$. Then, by
Lemma~\ref{waedfdrtdtgfz}, \begin{gather*} \widetilde{\cal
E}^{\mathrm B}_{\mu,n}({\cal W}^{(1)},{\cal W}^{(2)})= {\cal
E}^{\mathrm
B}_{\mu,n}(W^{(1)},W^{(2)})\\=\int_{\Gamma_X}\mu(d\gamma)\,\sum_{x\in\gamma}
\langle \nabla^X_x W^{(1)}(\gamma),\nabla^X_x
W^{(2)}(\gamma)\rangle_{T_xX\otimes
\wedge^n(T_\gamma\Gamma_X)}\notag\\= \sum_{k=1} ^n
\int_{\Gamma_X}\mu(d\gamma)\int_{X^k}{:}\,\gamma ^{\otimes
k}\,{:}(dx_1,\dots ,dx_k)\, \sum_{x\in\gamma} \langle \nabla^X_x
W^{(1)}_k(\gamma,x_1,\dots,x_k),\notag\\ \nabla^X_x
W^{(2)}_k(\gamma,x_1,\dots,x_k)\rangle_{T_xX\otimes {\Bbb
T}^{(n)}_{\{x_1,\dots,x_k\}}X^k}\notag\\ =\sum_{k=1}^n
\int_{\Gamma\times X^k}
\mu^{(k)}(d\gamma,dx_1,\dots,dx_k)\sum_{x\in\gamma\cup\{x_1,\dots,x_k\}}
\langle \nabla^X_x {\cal W}^{(1)}(\gamma,x_1,\dots,x_k),\notag\\
\nabla^X_x {\cal
W}^{(2)}(\gamma,x_1,\dots,x_k)\rangle_{T_xX\otimes {\Bbb
T}^{(n)}_{\{x_1,\dots,x_k\}}X^k},\end{gather*} which is equal to
the right hand side of \eqref{5443}.\quad
$\blacksquare$\vspace{2mm}

We will now apply Proposition~\ref{trdrtsers} to prove the
vanishing of square-integrable Bochner harmonic forms.

\setcounter{theorem}{8}
\begin{theorem}\label{rtsreaweraewra}
Let the conditions of Theorem~\rom{\ref{th4.1}} be satisfied\rom,
let \begin{equation}\label{rsresare}
\sigma^{(k)}(\gamma,X^k)=\infty\qquad \text{\rom{for $\mu$-a.e.\
$\gamma\in\Gamma_X$,}}\ k\in\N,\end{equation} and let one of the
two following conditions hold\rom:

\begin{description}

\item[\rom{(i)}] for $\mu$-a\rom.e\rom.\ $\gamma\in\Gamma_X$\rom,
$\rho(\gamma,\cdot)$ is continuous and positive on $X$\rom;

\item[\rom{(ii)}] $d\ge 2$ and for $\mu$-a\rom.e\rom.\
$\gamma\in\Gamma_X$\rom, $\rho(\gamma,\cdot)$ is continuous and
positive on $X\setminus\gamma$\rom.

\end{description}

 Then\rom, for each $n\in\N$\rom, $\operatorname{Ker} {\bf
H}_{\mu,n}^{\mathrm B}=\{0\}$\rom.

\end{theorem}

\noindent{\it Proof}. We will prove the theorem in the case of
(ii), the case (i) being completely similar and simpler.

First, we note that we can suppose that, for all
$\gamma\in\Gamma_X$, $\rho(\gamma,\cdot)$ is continuous and
positive on $X\setminus\gamma$. It suffices to show that ${\cal
E}^{\mathrm B}_{\mu,n}(W)=0$, $W\in D({\cal E}^{\mathrm
B}_{\mu,n})$ $\Rightarrow$ $W=0$, or equivalently,
$\widetilde{\cal E}^{\mathrm B}_{\mu,n}({\cal W})=0$, ${\cal W}\in
D(\widetilde{\cal E}^{\mathrm B}_{\mu,n})$ $\Rightarrow$ ${\cal
W}=0$. Here and below, for a bilinear form $E$ we set
$E(W){:=}E(W,W)$ for $W\in D(E)$.

Let us consider the following bilinear form on the Hilbert space
\eqref{tensor-n}: \begin{gather*} {\cal U}_n({\cal W}^{(1)},{\cal
W}^{(2)}){:=}\sum_{k=1}^n  {\cal U}_{k,n}({\cal W}^{(1)},{\cal
W}^{(2)}),\\ {\cal U}_{k,n}({\cal W}^{(1)},{\cal
W}^{(2)}){:=}\int_{\Gamma_X\times X^k}
\mu^{(k)}(d\gamma,dx_1,\dots,dx_k)\times\\ \times \langle
\nabla^{X^k}_{(x_1,\dots,x_k)}{\cal
W}^{(1)}(\gamma,x_1,\dots,x_k),
\nabla^{X^k}_{(x_1,\dots,x_k)}{\cal
W}^{(2)}(\gamma,x_1,\dots,x_k)\rangle_{T_{(x_1,\dots,x_k)}X^k\otimes
{\Bbb T}^{(n)}_{\{x_1,\dots,x_k\}}X^k},\\ {\cal W}^{(1)},{\cal
W}^{(2)}\in I^n({\cal D}\Omega^n)\end{gather*}  From the existence
of the generator of ${\cal U}_n$ defined on $I^n({\cal
D}\Omega^n)$, it follows that ${\cal U}_n$ is closable and let
$({\cal U}_n,D({\cal U}_n))$ denote its closure. By
Proposition~\ref{trdrtsers}, $D(\widetilde{\cal
E}_{\mu,n}^{\mathrm B})\subset D({\cal U}_n)$ and $\widetilde{\cal
E}_{\mu,n}^{\mathrm B}({\cal W})\ge {\cal U}_n({\cal W})$ for all
${\cal W}\in D(\widetilde{\cal E}_{\mu,n}^{\mathrm B})$.  Furthermore, it follows from the definition of ${\cal
U}_n$ that $D({\cal U}_n)=\bigoplus_{k=1}^n D({\cal U}_{k,n})$ and
${\cal U}_n=\sum_{k=1}^n {\cal U}_{k,n}$, where for each
$k=1,\dots,n$ $({\cal U}_{k,n},D({\cal U}_{k,n}))$ is a closed
form on $L^2_\Psi(\Gamma_X\times X^k\to
\wedge^n(TX^k);\mu^{(k)}){=:}H_{k,n}$. Hence, it suffices to show
that ${\cal U}_{k,n}({\cal W})=0$, ${\cal W}\in D({\cal U}_{k,n})$
$\Rightarrow$ ${\cal W}=0$.

For ${\cal W}\in I^n({\cal D}\Omega^n)\cap
H_{k,n}{=:}\Omega_{k,n}$, we define $$S({\cal
W})(\gamma,x_1,\dots,x_k){:=}
\|\nabla^{X^k}_{(x_1,\dots,x_k)}{\cal
W}(\gamma,x_1,\dots,x_k)\|^2$$ (here and below we omit the
notation of the space in the norm if this space is clear from the
context). Let $\{{\cal W}^{(n)}\}\subset\Omega_{k,n}$ and let
${\cal W}^{(n)}\to {\cal W}$ as $n\to\infty$ in the norm
$$\|\cdot\|_{D({\cal
U}_{k,n})}{:=}\big(\|\cdot\|_{H_{k,n}}^2+{\cal
U}_{k,n}(\cdot)\big)^{1/2}.$$ Using the inequality $$ \big(
S({\cal W}^{(n)})^{1/2}-S({\cal W}^{(m)})^{1/2}\big)^{2}\le
S({\cal W}^{(n)}-{\cal W}^{(m)}),$$ we conclude that $\{S({\cal W
}^{(n)})\}$ is a Cauchy sequence  in the norm of
$L^1(\Gamma_X\times X^k;\mu^{(k)} )$. Let $S({\cal W})$ denote its
limit. Then, using the definition of $\mu^{(k)}$,  we have
\begin{equation}\label{kjuhu} {\cal U}_{k,n}({\cal W})=
\int_{\Gamma_X}\mu(d\gamma)\int_{X^k} m(dx_1)\dotsm m(dx_k)\,
\rho^{(k)}(\gamma,x_1,\dots,x_k) S({\cal
W})(\gamma,x_1,\dots,x_k),\end{equation} where $$
\rho^{(k)}(\gamma,x_1,\dots,x_k){:=}\rho(\gamma,x_1)\rho(\gamma+\eps_{x_1},x_2)\dotsm
\rho(\gamma+\eps_{x_1}+\dots+\eps_{x_{k-1}},x_k).$$

Suppose now that ${\cal U}_{k,n}({\cal W})=0$.  Then, by  (ii), it
follows from \eqref{kjuhu} that, for $\mu$-a.e.\
$\gamma\in\Gamma_X$,
\begin{equation}\label{serasrwas} S({\cal W})(\gamma,\cdot)=0\quad
\text{$m^{\otimes k}$-a.e.\ on $X^k$}. \end{equation}

Let us fix $\gamma\in\Gamma_X$ such that \eqref{serasrwas} holds
and let ${\cal O}$ be an open ball in $X^k$ such that
\begin{equation}\label{reserws}\overline{\cal O}\subset {\cal
X}_{k,\gamma}{:=}\widetilde X^k\cap
(X\setminus\gamma)^k.\end{equation} Since
$\rho^{(k)}(\gamma,\cdot)$ is positive and continuous  on
$\overline{\cal O}$, $$ 0<c_1\le \rho^{(k)}(\gamma,\cdot)\le
c_2<\infty\quad \text{on }{\cal O}, $$ and so  $L^p$-convergence
on $\cal O$ with respect to the measure
$\sigma^{(k)}(\gamma,dx_1,\dots,dx_k)$ is equivalent to the same
convergence with respect to the measure $m^{\otimes k}$.

Let $W_2^1({\cal O})$ denote the Sobolev space consisting of all
functions $f\in L^2({\cal O};m^{\otimes k})$ which are weakly
differentiable and whose weak gradient $\nabla^{X^k} f\in
L^2({\cal O}\to T{\cal O};m^{\otimes k})$.

\setcounter{lemma}{9}
\begin{lemma}\label{easwa7644} We have $\|{\cal W}(\gamma,\cdot)\|\in W_2^1({\cal
O})$ and $\nabla^X \|{\cal W}(\gamma,\cdot)\|=0$ $m^{\otimes
k}$-a.e.\ on ${\cal O}$\rom. \end{lemma}

\noindent {\it Proof}. Let us consider the classical pre-Dirichlet
form on $L^2({\cal O};m^{\otimes k})$: $${\cal
E}(f^{(1)},f^{(2)})=\int_{{\cal O }}\langle \nabla^{X^k}
f^{(1)}(x_1,\dots,x_k),\nabla^{X^k} f^{(2)}(x_1,\dots
x_k)\rangle_{T_{(x_1,\dots,x_k)}X^k}\, m(dx_1)\dotsm m(dx_k),$$
where $f^{(1)}, f^{(2)}\in D({\cal E }){:=}C^1(\overline{\cal
O})$. As well known, this pre-Dirichlet form is closable and let
$({\cal E},D({\cal E}))$ denote its closure. Then, $D({\cal
E})=W_2^{1}({\cal O})$ and $$ {\cal E}(f^{(1)},f^{(2)})=\int_{\cal
O}{\cal S}(f^{(1)},f^ {(2)})(x_1,\dots,x_k)\, m(dx_1)\dotsm
m(dx_k),\qquad f^{(1)},f^{(2)}\in D({\cal E}),$$ where $$ {\cal
S}(f^{(1)},f^ {(2)})(x_1,\dots,x_k)= \langle \nabla^{X^k}
f^{(1)}(x_1,\dots,x_k),\nabla^{X^k}
f^{(2)}(x_1,\dots,x_k)\rangle_{T_{(x_1,\dots,x_k)}X^ k},$$ the
gradient $\nabla^{X^k}$ being understood in the weak sense.

Hence, taking notice of \eqref{serasrwas}, to prove  this lemma,
it suffices to show that the following claim is true: Let
$\omega:{\cal O}\to \wedge^n(T{\cal O})$  be a limit of a sequence
$\{\omega_n\}$ of smooth $n$-forms on $\overline{\cal O}$ with
respect to the norm $ \big( \|\cdot\|_{L^2({\cal O}\to \wedge^n
(T{\cal O });m^{\otimes k})}^2+ {\cal G}(\cdot)\big)^{1/2}$, where
$$ {\cal G}(u){:=}\int_{\cal
O}\|\nabla^{X^k}u(x_1,\dots,x_k)\|^2\, m(dx_1)\dotsm m(dx_k )$$
for a smooth form $u$. Then, $\|\omega\|\in D({\cal E})$ and
\begin{equation}\label{623}{\cal S}(\|\omega\|)\le S(\omega)\qquad \text{$m^{\otimes k}$-a.e.\ on
${\cal O}$}.\end{equation} Here, $S(\omega)(x_1,\dots,x_k)$ is
constructed analogously to the $S(W)(\gamma,x_1,\dots,x_k)$ above.

The proof of this claim is essentially the same as the proof of
the fact that, for each $f\in D({\cal E})$, $|f|\in D({\cal E})$
and ${\cal S}(|f|)\le {\cal S}(f)$ $m^{\otimes k}$-a.e., which is
why we limit ourselves to only outline it. So, first one shows by
approximation that, for each fixed $\epsilon>0$,
$\sqrt{\langle\omega,\omega\rangle+\epsilon}\in D({\cal E})$, and
moreover, for any fixed $\epsilon,\epsilon'>0$,
\begin{multline}\label{ztzdrs} {\cal S}
(\sqrt{\langle\omega,\omega\rangle+\epsilon}-\sqrt{\langle\omega,\omega\rangle+\epsilon'}\,)(x_1,\dots,x_k)\\
\le {\cal
S}(\omega)(x_1,\dots,x_k)\bigg\|\frac{\omega(x_1,\dots,x_k)}{\sqrt{\langle\omega,\omega\rangle+\epsilon}}-
\frac{\omega(x_1,\dots,x_k)}{\sqrt{\langle\omega,\omega\rangle+\epsilon'}}\bigg\|^2\quad
\text{$m^{\otimes k}$-a.e.\ }(x_1,\dots,x_k)\in{\cal O}.
\end{multline}
Second, one sets $\epsilon_n\downarrow 0$ and shows using
\eqref{ztzdrs} that
$\big\{\sqrt{\langle\omega,\omega\rangle+\epsilon_n}\,\big\}$ is a
Cauchy sequence with respect to the norm $\big(
\|\cdot\|_{L^2({\cal O};m^{\otimes k})}^2+{\cal
E}(\cdot)\big)^{1/2}$. The estimate \eqref{623} then trivially
follows. Thus, the lemma is proved.\quad $\blacksquare$

By Lemma~\ref{easwa7644}, it follows that $\|{\cal
W}(\gamma,\cdot)\|=\operatorname{const}$ $m^{\otimes k}$-a.e.\ on
$\cal O$. Since $d\ge2$, the set ${\cal X}_{k,\gamma}$ defined in
\eqref{reserws} is open and connected, and therefore it can be
covered by a countable number of open balls $\{{\cal O}_n\}$
satisfying $\overline{\cal O}_n\subset{\cal X}_{k,\gamma}$.
Therefore, $\|{\cal W}(\gamma,\cdot)\|=\operatorname{const}$
$m^{\otimes k}$-a.e.\ on ${\cal X}_{k,\gamma}$, and hence
$m^{\otimes k}$-a.e.\ on $X^k$. Finally, by \eqref{rsresare},
$\|W\|=0$ $\mu\otimes m^{\otimes k}$-a.e.\ on $\Gamma_X\times
X^k$. Thus, the theorem is proved.\quad $\blacksquare$

\subsection{deRham Laplacian on forms }

Let ${\cal E}\Omega^n$ denote the subset of ${\cal F}\Omega^n$
consisting of all forms $W\in{\cal F}\Omega^n$ such that all
derivatives of $W$ are polynomially bounded, that is, for each
$k\in\N$ there exist $\varphi\in{\cal D}$, $\varphi\ge0$, and
$l\in\N$ (depending on $W$) such that
\begin{equation}\label{0989454}
\|(\nabla^\Gamma)^{(k)}W(\gamma)\|^2_{(T_\gamma\Gamma_X)^{\otimes
k}\otimes\wedge^n(T_\gamma\Gamma_X)} \le \langle
\varphi,\gamma\rangle^l\qquad\text{for all
}\gamma\in\Gamma_X,\end{equation} and additionally, for each fixed
$\gamma\in\Gamma_X$ and $r\in\N$, the mapping $$
(X\setminus\gamma)^r\cap \widetilde X{}^r\ni(x_1,\dots,x_r)\mapsto
W(\gamma+\eps_{x_1}+\dots+\eps_{x_r})\in\wedge^n(T_\gamma\Gamma_X\oplus
T_{x_1}X\oplus\dots\oplus T_{x_r}X )$$ extends to a smooth form
$$X^r\ni(x_1,\dots,x_r) \mapsto \omega(x_1,\dots,x_r)\in
\wedge^n(T_\gamma\Gamma_X\oplus T_{x_1}X\oplus\dots\oplus T_{x_r}X
).$$ (Notice that the locality of a form, together with the above
condition of extension, will automatically imply the infinitely
differentiability of the form.)

As easily seen, ${\cal D}\Omega^n$ is a subset of ${\cal
E}\Omega^n$, and so we get the following  chain of inclusions
$${\cal D}\Omega^n\subset{\cal E}\Omega^n\subset{\cal
F}\Omega^n.$$

 We define linear operators
\begin{equation}\label{awertudrz}
{\bf d}_n\colon {\cal E}\Omega ^n\to {\cal E}\Omega ^{n+1},\qquad n\in {\Bbb %
Z}_+,\ {\cal E}\Omega ^0:={\cal F}C_{\mathrm b}^\infty({\cal D
},\Gamma_X),
\end{equation}
by
\begin{equation}\label{awgfse}
({\bf d}_nW)(\gamma ):=(n+1)^{1/2}\,\operatorname{AS}_{n+1}(\nabla
^\Gamma W(\gamma )),
\end{equation}
where \begin{equation}\label{uzui}\operatorname{AS}_{n+1}\colon
(T_\gamma \Gamma _X)^{\otimes (n+1)}\to \wedge ^{n+1}(T_\gamma
\Gamma _X)\end{equation} is the antisymmetrization operator. (We
notice that the polynomial boundedness of the form ${\bf d}_nW$
and its derivatives follows from the corresponding boundedness of
$\nabla^\Gamma W$ and the fact that the norm of the operator
\eqref{uzui} for each $\gamma\in\Gamma_X$ is equal to one).

Let us now consider ${\bf d}_n$ as an operator acting from the
space $L_\mu ^2\Omega ^n$ into $L_\mu ^2\Omega ^{n+1}$. (We remark
that, by the proof of Lemma~\ref{uhizhgzt}, ${\bf d}_nW=0$
$\mu$-a.e.\ for $W\in {\cal E}\Omega^n$ such that $W=0$
$\mu$-a.e.) We denote by ${\bf d}^*_{n}$ the adjoint operator of
${\bf d}_n$.

\setcounter{proposition}{10}
\begin{proposition}\label{guzftdrt} Let \eqref{awrzwerjidrt}
hold\rom. Then\rom,
 ${\bf d}^*_{n}$ is a densely defined operator from
$L^2_\mu\Omega^{n+1}$ into $ L_\mu^2\Omega^n$ with domain
containing ${\cal E}\Omega ^{n+1}$\rom.
\end{proposition}

\noindent {\it Proof}. It follows from \eqref{awgfse} and the
definition of $\nabla^\Gamma$ that, for any $W\in{\cal E}\Omega^n$
and $\gamma\in\Gamma_X$, \begin{equation}\label{23903}({\bf d}_n
W) (\gamma)=\sum_{x\in\gamma}({\bf
d}_{x,n}W)(\gamma),\end{equation} where
\begin{equation}\label{qwk86}({\bf d}_{x,n}W)(\gamma):=(n+1)^{1/2}
\operatorname{AS}_{n+1}(\nabla^X_x W(\gamma)).\end{equation}

Let $\gamma \in \Gamma _X$ and $x\in \gamma $ be fixed. Let $%
C^\infty ({\cal O}_{\gamma ,x}\rightarrow \wedge ^n(T_\gamma
\Gamma _X))$ denote the space of all smooth sections of the
Hilbert bundle (\ref{bund1}). We define an operator
\[
d_{x,n}^X:C^\infty ({\cal O}_{\gamma ,x}\to \wedge ^n(T_\gamma
\Gamma _X))\to C^\infty ({\cal O}_{\gamma ,x}\rightarrow \wedge
^{n+1}(T_\gamma \Gamma _X))
\]
whose action, in local coordinates  on the manifold $X$, is given
by
\begin{equation}\label{aweese}
d^X_{x,n}\,\phi (y)\,h_1\wedge \dots \wedge
h_n=(n+1)^{1/2}\,\nabla ^X\phi (y)\wedge h_1\wedge \dots \wedge
h_n,
\end{equation}
$\phi \in C^\infty ({\cal O}_{\gamma ,x}\to\R)$, $h_k\in
T_{x_k}X$, $x_k\in \gamma$, $k=1,\dots,n$. It follows  from
\eqref{qwk86} and \eqref{aweese} that
\begin{equation}\label{798cdjhgztf} ({\bf
d}_{x,n}W)(\gamma)=d^X_{x,n}W_x(\gamma,x).\end{equation}

Next, let $\Omega({\cal O}_{\gamma ,x}\rightarrow \wedge
^n(T_\gamma \Gamma _X))$ denote the space of all  sections of the
Hilbert bundle (\ref{bund1}). We define an operator
\[
\delta^X _{x,n}:C^\infty ({\cal O}_{\gamma ,x}\rightarrow \wedge
^{n+1}(T_\gamma \Gamma _X))\rightarrow \Omega({\cal O}_{\gamma
,x}\rightarrow \wedge ^n(T_\gamma \Gamma _X))
\]
setting
\begin{multline}
\delta ^X_{x,n}\,\phi (y)\,h_1\wedge \dots \wedge
h_{n+1}{:=}-(n+1)^{-1/2}\,\sum_{i=1}^{n+1}(-1)^{i-1}\varepsilon_{x,x_i}
\big[ \langle \nabla ^X\phi (y),h_i\rangle _x \\
\text{}+\phi(y)\langle B_\mu (\gamma,y),h_i\rangle_x\big]
h_1\wedge \dots \wedge \check{h}_i\wedge \dots \wedge
h_{n+1},\label{awetz}
\end{multline}
where $\phi \in C^\infty ({\cal O}_{\gamma ,x}\to\R)$, $h_i\in
T_{x_i}X$, $x_i\in \gamma$, $i=1,\dots,n+1$,
\[
\varepsilon_{x,x_i}{:=}\begin{cases}1,& x=x_i,\\ 0,&
\text{otherwise},\end{cases}
\]
and $\check{h}_i$ denotes the absence of $h_i$. We now set for $W\in {\cal %
E}\Omega ^{n+1}$%
\begin{equation}\label{oipdrt}
({\pmb \delta }_{x,n}W)(\gamma){:=}{\delta}^X_{x,n}W_x(\gamma,x)
\end{equation}
and \begin{equation}\label{tz7634} ({\pmb \delta
}_{n}W)(\gamma){:=}\sum_{x\in\gamma}({\pmb \delta
}_{x,n}W)(\gamma)\end{equation} (notice that the sum on the right
hand side of \eqref{tz7634} is actually finite).

Let us show that, for any $W\in {\cal E}\Omega^{n+1}$, we have
${\pmb \delta}_{n}W\in L^2_\mu\Omega^n$, where
$L^2_\mu\Omega^0{:=}L^2(\Gamma_X;\mu)$. We choose any $\Lambda\in
{\cal O}_c(X)$ such that $W(\gamma)=W(\gamma_{\Lambda'})$ for some
compact  $\Lambda'\subset\Lambda$. Then, by \eqref{tz7634},
\begin{align}\label{wahges78}&\int_{\Gamma_X}\|({\pmb\delta}_{n}W)(\gamma)\|^2_{\wedge^n(T_\gamma\Gamma_X)}
\,\mu(d\gamma)\\&\qquad =\int_{\Gamma_X}\left\|
\sum_{x\in\gamma_\Lambda}({\pmb\delta}_{x,n}W)(\gamma)\right\|^2_{\wedge^n(T_\gamma\Gamma_X)}
\,\mu(d\gamma)\notag\\ &\qquad
\le\int_{\Gamma_X}\left(\sum_{x\in\gamma_\Lambda}
\|({\pmb\delta}_{x,n}W)(\gamma)\|_{\wedge^n(T_\gamma\Gamma_X)}\right)^2\,\mu(d\gamma).
\notag
\end{align} Using \eqref{0989454},  \eqref{awetz},  and \eqref{oipdrt},
it is not hard to show that there exist $\varphi\in C_0(X)$,
$\varphi\ge0$, and $k\in\N$ (independent of $\gamma$ and $x$) such
that
\begin{equation}\label{875243} \|({\pmb
\delta}_{x,n}W)(\gamma)\|_{\wedge^n(T_\gamma\Gamma_X)}\le \langle
\varphi,\gamma\rangle^k+(n+1)^{1/2}|B_\mu(\gamma,x)|_x\|W(\gamma)\|
_{\wedge^{n+1}(T_\gamma\Gamma_X)}.\end{equation} Analogously to
the proof of \eqref{tfujipp}, we get from \eqref{weghtse},
\eqref{edrdt} \eqref{awrzwerjidrt}, \eqref{wahges78}, and
\eqref{875243} that ${\pmb \delta}_{n}W\in L^2_\mu\Omega^n$.

Let $W^{(1)},W^{(2)}\in{\cal E}\Omega^n$ and let $\Lambda\in{\cal
O} _c(X)$ be such that, for some compact $\Lambda'\subset\Lambda$
$W^{(i)}(\gamma)=W^{(i)}(\gamma_{\Lambda'})$, $i=1,2$, for all
$\gamma\in\Gamma_X$. Then, by \eqref{qwgfhjzjg},  \eqref{23903},
\eqref{aweese}, and \eqref{798cdjhgztf},  we get using the
notations of Section~\ref{ewew086}
\begin{gather}
\int_{\Gamma_X}\langle {\bf
d}_nW^{(1)}(\gamma),W^{(2)}(\gamma)\rangle_{\wedge^{n+1}(T_\gamma\Gamma_X)}\,\mu(d\gamma)
\label{ewewee}\\
=\int_{\Gamma_X}\mu(d\gamma)\int_\Lambda\gamma(dx)\, \langle ({\bf
d}_{x,n}W^{(1)})(\gamma),W^{(2)}(\gamma)\rangle_{\wedge^{n+1}(T_\gamma\Gamma_X})\notag\\
=\int_{\Gamma_X}\mu(d\gamma)\int_\Lambda\sigma(\gamma,dx)\,
\langle ({\bf d}_{x,n}W^{(1)})(\gamma+\eps_x),
W^{(2)}(\gamma+\eps_x)\rangle_{\wedge^{n+1}(T_{\gamma+\eps_x}\Gamma_X})\notag\\
=\int_{\Gamma_X}\mu(d\gamma)\int_\Lambda\sigma(\gamma,dx)\,
\sum_{k=1}^n \sum_{\begin{gathered} \scriptstyle
\{x_1,\dots,x_k\}\subset\gamma\cup\{x\}
\\ \scriptstyle x\in\{x_1,\dots,x_k\}
\end{gathered}}
 \la ({\bf
d}_{x,n}W^{(1)})_k(\gamma+\eps_x,x_1,\dots,x_k),\notag\\
W^{(2)}_k(\gamma+\eps_x,x_1,\dots,x_k)\ra_
{\wedge^{n+1}(T_{\gamma+\eps_x}\Gamma_X)}\notag\\
=\int_{\Gamma_X}\mu(d\gamma) \sum_{k=1}^n
\sum_{\{x_1,\dots,x_{k-1}\}\subset\gamma}\int_\Lambda\sigma(\gamma,dx)\,
\la ({\bf
d}_{x,n}W^{(1)})_k(\gamma+\eps_x,x,x_1,\dots,x_{k-1}),\notag\\
W^{(2)}_k(\gamma+\eps_x,x,x_1,\dots,x_{k -1})\ra_
{\wedge^{n+1}(T_{\gamma+\eps_x}\Gamma_X)}.\notag
\end{gather}
It follows from the definition of ${\cal E}\Omega^n$ that, for a
fixed $\gamma\in\Gamma_X$ and
$\{x_1,\dots,x_{k-1}\}\subset\gamma$,
$W_k^{(2)}(\gamma,\cdot,x_1,\dots,x_{k-1})$ extends to a smooth
form \begin{multline*}X\ni x\mapsto
W^{(2)}_k(\gamma,x,x_1,\dots,x_{k-1})\in\\ \in \bigoplus_{\begin{gathered}{%
\scriptstyle{ 1\le l_1,\dots,l_k\le d}} \\ \scriptstyle
l_1+\dots+l_k=n
\end{gathered}}(T_xX)^{\wedge l_1}\wedge (T_{x_1}X)^{\wedge l_2}\wedge
\dots \wedge (T_{x_{k-1}}X)^{\wedge l_k}\subset \big(T_xX\oplus
T_{x_1}X\oplus T_{x_{k-1}}X \big)^{\wedge n}.\end{multline*} Since
$W^{(1)}(\gamma+\eps_{\bullet})$ also extends to a smooth form on
$X$, we can carry out an integration by parts in the $x$ variable
in   \eqref{ewewee}. Thus, by using \eqref{awetz}--\eqref{tz7634}
and \eqref{qwgfhjzjg}, we continue \eqref{ewewee} as follows:
\begin{gather}= \int_{\Gamma_X}\mu(d\gamma)\sum_{k=1}^n
\sum_{\{x_1,\dots,x_{k-1}\}\subset\gamma}\int_\Lambda
\sigma(\gamma,dx)\la W^{(1)}(\gamma+\eps_x),\notag\\
\delta^X_{x,n}W^{(2)}_k
(\gamma+\eps_x,x,x_1,\dots,x_{k-1})\rangle_{\wedge^{n}(T_{\gamma+\eps_x}\Gamma_X)}\notag\\
=
\int_{\Gamma_X}\mu(d\gamma)\int_\Lambda\sigma(\gamma,dx)\, \langle
W^{(1)} (\gamma+\eps_x),\delta^X_{x,n} W_x^{(2)}(\gamma+\eps_x,x)
\rangle_{\wedge^{n}(T_{\gamma+\eps_x}\Gamma_X)}\notag\\
=
\int_{\Gamma_X}\mu(d\gamma)\int_\Lambda\sigma(\gamma,dx)\, \langle
W^{(1)} (\gamma+\eps_x),({\pmb \delta}_{x,n}
W^{(2)})(\gamma+\eps_x)
\rangle_{\wedge^{n}(T_{\gamma+\eps_x}\Gamma_X)}\notag\\
=\int_{\Gamma_X}\mu(d\gamma)\int_\Lambda \gamma(dx)\, \langle
W^{(1)}(\gamma),({\pmb\delta}_{x,n}W^{(2)})(\gamma)
\rangle_{\wedge^n(T_\gamma\Gamma_X)}\notag\\
=\int_{\Gamma_X}\langle
W^{(1)}(\gamma),({\pmb\delta}_{n}W^{(2)})(\gamma)
\rangle_{\wedge^n(T_\gamma\Gamma_X)} \,\mu(d\gamma).
\notag
\end{gather}

Hence, ${\cal F}\Omega^{n+1}\subset D({\bf d}^*_{\mu,n})$ and
${\bf d}^*_{n}\restriction{\cal E}\Omega^{n+1}
={\pmb\delta}_{\mu,n}$.\quad $\blacksquare$ \vspace{2mm}

\setcounter{corollary}{11}
\begin{corollary} The operator ${\bf d}_n: L^2_\mu\Omega^n\to L^2_\mu\Omega^{n+1}$ is
closable\rom.\end{corollary}

We denote  by $\bar{\bf d}_n$ the closure of ${\bf d}_n$.
The space $%
Z^n:=\operatorname{Ker}\bar{\bf d}_n$ is then a closed subspace of
$L_\pi ^2\Omega ^n$. Let  $B^n$ denote the closure in $L_\pi
^2\Omega ^n$ of the subspace $\operatorname{Im}{\bf d}_{n-1}$ (of
course, $B^n=$the closure of $\operatorname{Im}\bar{\bf
d}_{n-1}$).

We obviously have
$
{\bf d}_n{\bf d}_{n-1}=0 $, which implies $$ \operatorname{Im}{\bf
d}_{n-1}\subset\operatorname{Ker}{\bf d}_n\subset Z^n.$$ Hence
$B^n\subset Z^n$ and \begin{equation}\label{pserdo}
 \bar{\bf d}_n\bar{\bf
d}_{n-1}=0.\end{equation}

Thus, we have the infinite complex
\[
 \cdots\stackrel{{\bf d}_{n-1}}{\longrightarrow }{\cal E}\Omega ^n\stackrel{{\bf d}_n%
}{\longrightarrow }{\cal E}\Omega ^{n+1}\stackrel{{\bf d}_{n+1}}{\longrightarrow }%
\cdots\, ,
\]
 and the associated Hilbert complex
\begin{equation}
\cdots\stackrel{\bar{\bf d}_{n-1}}{\longrightarrow }L_\pi ^2\Omega ^n\stackrel{%
\bar{\bf d}_n}{\longrightarrow }L_\pi ^2\Omega ^{n+1}\stackrel{\bar{\bf d}%
_{n+1}}{\longrightarrow }\cdots\, .  \label{complex}
\end{equation}
 We set in a standard way
\[
{\cal H}_\mu ^n=Z^n/B^n,\qquad n\in\N.\]

For $n\in {\Bbb N}$, we define a bilinear form ${\cal E}_{\mu,n}
^{\mathrm R}$ on $L^2_\mu\Omega^n$ by
\begin{multline}
{\cal E}_{\mu, n}
^{\mathrm R}(W^{(1)},W^{(2)}):=\int_{\Gamma _X}\big[ \langle \bar{\bf d}%
_nW^{(1)}(\gamma ),\bar{\bf d}_nW^{(2)}(\gamma )\rangle _{\wedge
^{n+1}(T_\gamma \Gamma _X)}  \label{lklk} \\
\text{{}}+\langle {\bf d}_{n-1}^{*}W^{(1)}(\gamma ),{\bf d}%
_{n-1}^{*}W^ {(2)}(\gamma )\rangle _{\wedge ^{n-1}(T_\gamma \Gamma _X)}\big]%
\,\mu (d\gamma ),
\end{multline}
where $W^{(1)},W^{(2)}\in D({\cal E}_{\mu,n} ^{\mathrm
R}){:=}D(\bar{\bf d}_n)\cap D({\bf d}^*_{n-1})$. This form is
evidently closed, and let \linebreak $({\bf H}_{\mu,n}^{\mathrm
R},D({\bf H}_{\mu,n}^{\mathrm R}))$ denote its generator. This
operator  will be called {\it the Hodge--deRham Laplacian of the
measure $\mu$}.

The following proposition reflects a quite standard fact in the
theory of $L^2$-cohomologies.

\setcounter{proposition}{12}
\begin{proposition}
The natural isomorphism between ${\cal H}^n_\mu$ and the
orthogonal complement of $B^n$ to $Z^n$ is the isomorphism of the
Hilbert spaces
\begin{equation}
{\cal H}_\mu ^n\simeq \operatorname{Ker}{\bf H}_{\mu,n}^{\mathrm
R}. \label{harm1}
\end{equation}
\end{proposition}

\noindent {\it Proof}. Using \cite[Proposition~A.1]{Ar0}, we
conclude from Proposition~\ref{guzftdrt} and formula
\eqref{pserdo} that
\begin{equation}
L_\mu ^2\Omega ^n= \operatorname{Ker}{\bf H}_{\mu,n}^{\mathrm
R}\oplus \overline{\operatorname{Im} {\bf
d}_{n-1}}%
\oplus \overline{\operatorname{Im}{\bf d}_n^{*}}  \label{hdr}
\end{equation}
(the weak Hodge--deRham decomposition). For the closed operator
$\bar {\bf d}_n$ we have the standard decomposition
\[
L_\mu ^2\Omega ^n=\operatorname{Ker}\bar {\bf d}_{n}\oplus
\overline{\operatorname{Im} {\bf d}_{n}^*},
\]
which together with (\ref{hdr}) implies the result. $\blacksquare
$\vspace{2mm}

We do not know {\it a priori} whether the domain $D({\bf
H}_{\mu,n}^{\mathrm R})$ contains ${\cal D}\Omega^n$, however the
following theorem gives a sufficient condition for this.

\setcounter{theorem}{13}
\begin{theorem}\label{dtdrtdtdsrt} Let us suppose that\begin{description}

\item[\rom{(i)}] for $\mu$-a\rom.e\rom.\ $\gamma\in\Gamma_X$\rom,
$\rho(\gamma,x)>0$ for all $x\in X\setminus\gamma$ and the
function $\rho(\gamma,\cdot)$ is continuous on $X$\rom;

\item[\rom{(ii)}] for $\mu$-a\rom.e\rom.\ $\gamma\in\Gamma_X$\rom,
 $\rho(\gamma,\cdot)$ is two times  differentiable on
$X\setminus\gamma$ and  $\nabla^X\rho(\gamma,\cdot)$ extends to a
 continuous form on $X$\rom;

\item[\rom{(iii)}] for $\mu\otimes m $-a\rom.e\rom.\
$(\gamma,x)\in\Gamma_X\times X$\rom,  $y\mapsto
\nabla^X_x\rho(\gamma+\eps_y,x)\in T_xX$ is  differentiable on
$X\setminus(\gamma\cup\{x\})$  and
$$X\setminus(\gamma\cup\{x\})\ni y\mapsto
\frac{\rho(\gamma+\eps_x,y)}{\rho(\gamma+\eps_y,x)}\,\nabla^X_x\rho(\gamma
+\eps_y,x)\in T_xX$$ extends to a continuous mapping
 on $X$\rom;

\item[\rom{(iv)}] \eqref{awrzwerjidrt} holds\rom, and furthermore
\begin{equation}\label{awrzejirt}\forall\Lambda\in{\cal O}_c(X)\ \exists\eps>0:\quad
\int_{\Gamma_X}\left(\sum_{y\in\gamma}\sum_{x\in\gamma_\Lambda}\|\nabla^X_y
B_\mu(\gamma,x) \|_{T_yX\otimes
T_xX}\right)^{2+\eps}\,\mu(d\gamma)<\infty.\end{equation}

\end{description}
Then\rom, ${\cal D}\Omega^n\subset D({\bf H}_{\mu,n}^{\mathrm R})$
and  $${\bf H}^{\mathrm R}_{\mu,n}\restriction {\cal
D}\Omega^n={\bf d}_{n-1}{\bf d}_{n-1}^*+{\bf d}^*_{n}{\bf d}_n.$$
\end{theorem}

\noindent{\it Proof}. Since by Proposition~\ref{guzftdrt} ${\bf
d}_n{\cal D}\Omega^n\subset {\cal E}\Omega^{n+1}\subset D({\bf
d}^*_{n})$, to prove the theorem we have to show that ${\bf
d}^*_{n-1}{\cal D}\Omega^n\subset D({\bf d}_{n-1})$, that is, for
arbitrary $W^{(1)},W^{(2)}\in{\cal D }\Omega^n$, there exits $V\in
L^2_\mu\Omega^n$ such that
\begin{equation}\label{gftzdztdtrstrs} \int_{\Gamma_X}\langle {\bf
d}^*_{n-1}W^{(1)}(\gamma),{\bf d}_{n-1}^* W^{(2)}(\gamma)
\rangle_{\wedge^n(T_\gamma\Gamma_X)}\,\mu(d\gamma)=\int_{\Gamma_X}\langle
V(\gamma),
W^{(2)}(\gamma)\rangle_{\wedge^n(T_\gamma\Gamma_X)}\,\mu(d\gamma).\end{equation}

We choose any $\Lambda\in{\cal O}_c(X)$ such that, for some
compact $\Lambda'\subset\Lambda$,
$W^{(i)}(\gamma)=W^{(i)}(\gamma_{\Lambda'})$, $i=1,2$, for all
$\gamma\in\Gamma_X$. It follows from the proof of
Proposition~\ref{guzftdrt} and Lemma~\ref{waedfdrtdtgfz} that
\begin{gather} \label{9832432}\int_{\Gamma_X}\langle {\bf
d}^*_{n-1}W^{(1)}(\gamma),{\bf
d}^*_{n-1}W^{(2)}(\gamma)\rangle_{\wedge^{n-1}(T_\gamma\Gamma_X)}\,\mu(d\gamma)\\
=\int_{\Gamma_X}\sum_{x,y\in\gamma_\Lambda}
\langle{\pmb\delta}_{x,n}W^{(1)}(\gamma),{\pmb\delta}_{y,n}W^{(2)}(\gamma)\rangle
_{\wedge^{n-1}(T_\gamma\Gamma_X)}\,\mu(d\gamma)\notag\\=
\int_{\Gamma_X}\mu(d\gamma)\int_\Lambda
\sigma(\gamma,dx)\,\langle{\pmb\delta}_{x,n}W^{(1)}(\gamma+\eps_x)
,{\pmb\delta}_{x,n}W^{(2)}(\gamma+\eps_x)\rangle
_{\wedge^{n-1}(T_{\gamma+\eps_x}\Gamma_X)}\notag\\ \text{}+
\int_{\Gamma_X}\mu(d\gamma)\int_\Lambda
\sigma(\gamma,dx)\int_\Lambda \sigma(\gamma+\eps_x,dy)\,\langle
{\pmb\delta}_{x,n}W^{(1)}(\gamma+\eps_x+\eps_y),\notag\\
{\pmb\delta}_{y,n}W^{(2)}(\gamma+\eps_x+\eps_y)
\rangle_{\wedge^{n-1}(T_{\gamma+\eps_x+\eps_y} \Gamma_X)}
.\notag\end{gather}

Due to conditions (i)--(iii), we can see that, for $\mu$-a.e.\
$\gamma\in\Gamma_X$ and $x\in X\setminus\gamma$, $${\bf
d}_{x,n}{\pmb\delta}_{x,n} W^{(1)}(\gamma+\eps_x)\in\wedge
^n(T_{\gamma+\eps_x})$$ and for $\mu\otimes m$-a.e.\
$(\gamma,x)\in\Gamma_X\times X$ and $y\in
X\setminus(\gamma\cup\{x\})$, $${\bf
d}_{y,n}{\pmb\delta}_{x,n}W^{(1)}(\gamma+\eps_x+\eps_y)\in
\wedge^n(T_{\gamma+\eps_x+\eps_y}\Gamma_X), $$ using formulas
\eqref{aweese}, \eqref{798cdjhgztf},  for the definition of ${\bf
d}_{x,n}$, $x\in X$. Moreover, by virtue of (i) and (ii), the
integration by parts yields, for $\mu$-a.e.\ $\gamma\in\Gamma_X$
\begin{multline}\label{jeshg} \int_\Lambda
\sigma(\gamma,dx)\,\langle{\pmb\delta}_{x,n}W^{(1)}(\gamma+\eps_x)
,{\pmb\delta}_{x,n}W^{(2)}(\gamma+\eps_x)\rangle
_{\wedge^{n-1}(T_{\gamma+\eps_x}\Gamma_X)}\\ =
\int_{\Lambda}\sigma(\gamma,dx)\,\langle {\bf d}_{x,n}{\pmb
\delta}_{x,n}W^{(1)}(\gamma+\eps_x),W^{(2)}(\gamma+\eps_x)\rangle_{\wedge^n
(T_{\gamma+\eps_x}\Gamma_X)},\end{multline} and analogously, using
(i) and (iii), we get, for $\mu\otimes m$-a.e.\
$(\gamma,x)\in\Gamma_X\times X$,
\begin{multline}\label{jesh4554g}
\int_\Lambda \sigma(\gamma+\eps_x,dy)\,\langle
{\pmb\delta}_{x,n}W^{(1)}(\gamma+\eps_x+\eps_y),
{\pmb\delta}_{y,n}W^{(2)}(\gamma+\eps_x+\eps_y)
\rangle_{\wedge^{n-1}(T_{\gamma+\eps_x+\eps_y} \Gamma_X)}\\ =
\int_\Lambda \sigma(\gamma+\eps_x,dy)\,\langle {\bf
d}_{y,n}{\pmb\delta}_{x,n}W^{(1)}(\gamma+\eps_x+\eps_y),
W^{(2)}(\gamma+\eps_x+\eps_y)
\rangle_{\wedge^n(T_{\gamma+\eps_x+\eps_y}\Gamma_X)}.\end{multline}

Suppose that
\begin{equation}\label{o6564}\int_{\Gamma_X}\left(
\sum_{x,y\in\gamma}\|{\bf
d}_{y,n}{\pmb\delta}_{x,n}W^{(1)}(\gamma)\|
_{\wedge^n(T_\gamma\Gamma_X)}\right)^2\,\mu(d\gamma)<\infty,\end{equation}
so that \begin{equation*}V(\gamma){:=}
 \sum_{x,y\in\gamma}
{\bf d}_{y,n}{\pmb\delta}_{x,n}W^{(1)}(\gamma)\in
\wedge^n(T_\gamma\Gamma_X)\end{equation*} is well-defined for
$\mu$-a.a.\ $\gamma\in\Gamma_X$, and moreover $V\in
L^2_\mu\Omega^n$. Then, by Lemma~\ref{waedfdrtdtgfz} and
\eqref{jeshg}--\eqref{o6564}, we continue \eqref{9832432} as
follows:
\begin{gather}\label{zztfztd}
=\int_{\Gamma_X}\mu(d\gamma)\int_{\Lambda}\sigma(\gamma,dx)\,\langle
{\bf d}_{x,n}{\pmb
\delta}_{x,n}W^{(1)}(\gamma+\eps_x),W^{(2)}(\gamma+\eps_x)\rangle_{\wedge^n
(T_{\gamma+\eps_x}\Gamma_X)}\\ \text{}+
\int_{\Gamma_X}\mu(d\gamma)\,\int_\Lambda
\sigma(\gamma,dx)\int_\Lambda \sigma(\gamma+\eps_x,dy)\,\langle
{\bf
d}_{y,n}{\pmb\delta}_{x,n}W^{(1)}(\gamma+\eps_x+\eps_y),\notag\\
W^{(2)}(\gamma+\eps_x+\eps_y)
\rangle_{\wedge^n(T_{\gamma+\eps_x+\eps_y}\Gamma_X)}
=\int_{\Gamma_X}\langle
V_\Lambda(\gamma),W^{(2)}(\gamma)\rangle_{\wedge^n(T_\gamma\Gamma_X)}\,\mu(d\gamma),
\notag\end{gather} where
\begin{equation}\label{izrdrde}V_\Lambda(\gamma){:=}\sum_{x,y\in\gamma_\Lambda} {\bf
d}_{y,n}{\pmb\delta}_{x,n}W^{(1)}(\gamma).\end{equation}

Since $V_\Lambda(\gamma)\to V(\gamma)$ as $\Lambda\to X$ for
$\mu$-a.e.\ $\gamma\in\Gamma_X$, by the majorized convergence
theorem, we conclude from \eqref{9832432}, \eqref{zztfztd}, and
\eqref{izrdrde} that $$ \int_{\Gamma_X}\langle {\bf
d}^*_{n-1}W^{(1)}(\gamma),{\bf
d}^*_{n-1}W^{(2)}(\gamma)\rangle_{\wedge^{n-1}(T_\gamma\Gamma_X)}\,\mu(d\gamma)
=\int_{\Gamma_X}\langle
V(\gamma),W^{(2)}(\gamma)\rangle_{\wedge^n(T_\gamma\Gamma_X)}\,\mu(d\gamma).$$

Thus, it remains to show that \eqref{o6564} does indeed hold. Let
$\widetilde\Lambda\in{\cal O}_c(X)$ be such that, for some compact
$\Lambda'\subset\widetilde\Lambda$,
$W^{(1)}(\gamma)=W^{(1)}(\gamma_{\Lambda'})$ for all
$\gamma\in\Gamma_X$ ($\widetilde\Lambda$ being now independent of
$W^{(2)}$). Since ${\pmb \delta}_{x,n}W^{(1)}(\gamma)=0$ for all
$x\in\gamma_{\widetilde\Lambda{}^c}$, we get
\begin{multline}\label{rewsers}
\int_{\Gamma_X}\left( \sum_{x,y\in\gamma}\|{\bf
d}_{y,n}{\pmb\delta}_{x,n}W^{(1)}(\gamma)\|
_{\wedge^n(T_\gamma\Gamma_X)}\right)^2\,\mu(d\gamma)\\ \le
2\int_{\Gamma_X}\left(\sum_{x,y\in\gamma_{\widetilde\Lambda}}\|{\bf
d}_{y,n}{\pmb\delta}_{x,n}W^{(1)}(\gamma)\|
_{\wedge^n(T_\gamma\Gamma_X)}\right)^2\,\mu(d\gamma)\\ \text{}+
2\int_{\Gamma_X}\left(\sum_{y\in\gamma_{\widetilde\Lambda{}^{\mathrm
c }}}\sum_{x\in\gamma_{\widetilde\Lambda}}\|{\bf
d}_{y,n}{\pmb\delta}_{x,n}W^{(1)}(\gamma)\|
_{\wedge^n(T_\gamma\Gamma_X)}\right)^2\,\mu(d\gamma).\end{multline}

Analogously to \eqref{875243}, we conclude from
\eqref{aweese}--\eqref{oipdrt} the existence of $\varphi\in
C_0(X)$, $\varphi\ge0$, and $k\in\N$ (independent of $\gamma$,
$x$, and $y$) such that \begin{equation}\label{5454jg} \|{\bf
d}_{y,n}{\pmb\delta}_{x,n}W^{(1)}(\gamma)\|_{\wedge^n(T_\gamma\Gamma_X)}
\le
\langle\varphi,\gamma\rangle^k\big(
1+|B_\mu(\gamma,x)|_x+\|\nabla^X_y B_\mu(\gamma,x)\|_{T_yX\otimes
T_xX } \big) \end{equation} for
$x,y\in\gamma_{\widetilde\Lambda}$, and
\begin{equation}\label{seetrtr}
 \|{\bf
d}_{y,n}{\pmb\delta}_{x,n}W^{(1)}(\gamma)\|_{\wedge^n(T_\gamma\Gamma_X)}
\le
\langle\varphi,\gamma\rangle^k \|\nabla^X_y
B_\mu(\gamma,x)\|_{T_yX\otimes T_xX } ,\end{equation} for
$y\in\gamma_{\widetilde\Lambda{}^c}$ and
$x\in\gamma_{\widetilde\Lambda}$.

Thus, the finiteness of the right hand side of \eqref{rewsers} can
easily be deduced from \eqref{edrdt}, \eqref{awrzwerjidrt},
\eqref{awrzejirt}, \eqref{5454jg}, \eqref{seetrtr}, and the
Schwarz inequality. \quad $\blacksquare$

\setcounter{corollary}{14}
\begin{corollary} \label{pouoq213w}
Let the conditions of Theorem~\rom{\ref{dtdrtdtdsrt}} be
satisfied\rom. Then\rom, for each $W\in{\cal D}\Omega^n$ and
$\mu$-a\rom.e\rom.\ $\gamma\in\Gamma_X$
\begin{equation}\label{3w43276756}\sum_{x,y\in\gamma}\big(\| {\pmb\delta}_{x,n}{\bf d }
_{y,n}W(\gamma)\|_{\wedge^n(T_\gamma\Gamma_X)}+\|{\bf
d}_{y,n}{\pmb\delta}_{x,n}W(\gamma)\||_
{\wedge^n(T_\gamma\Gamma_X)}\big)<\infty,\end{equation} and the
action of the operator ${\bf H}_{\mu,n}^{\mathrm R}$ can be
represented in the form $${\bf H}^{\mathrm R}_{\mu,n}
W(\gamma)=\sum_{x,y\in\gamma} \big({\pmb\delta}_{x,n}{\bf d }
_{y,n}+{\bf d}_{y,n}{\pmb\delta}_{x,n}\big)W(\gamma),\qquad
\text{\rom{$\mu$-a.e.}\ $\gamma\in\Gamma_X$}. $$\end{corollary}

\subsection{Weitzenb\"ock formula}
In this section, we will derive a Weitzenb\"ock type formula,
which gives a relation between the Bochner Laplacian ${\bf
H}^{\mathrm B}_{\mu,n}$ and the deRham Laplacian ${\bf H}^{\mathrm
R}_{\mu,n}$. In what follows, we will suppose that  the conditions
of Theorem~\rom{\ref{dtdrtdtdsrt}} are satisfied.

For each $V(\gamma)\in T_\gamma\Gamma_X$, $\gamma\in\Gamma_X$, we
define an annihilation operator $$a_n(V(\gamma)):
\wedge^n(T_\gamma\Gamma_X)\to\wedge^{n-1}(T_\gamma\Gamma_X)  $$
and a creation operator $$a_n^*
(V(\gamma)):\wedge^{n-1}(T_\gamma\Gamma_X)\to\wedge^n(T_\gamma\Gamma_X)$$
as follows: \begin{align*} a_n(V(\gamma))W_n(\gamma)&{:=}\sqrt n\,
\langle V(\gamma),W_n(\gamma)\rangle_\gamma,\qquad W_n(\gamma)\in
\wedge^n(T_\gamma\Gamma_X),\\
a^*_n(V(\gamma))W_{n-1}(\gamma)&{:=}\sqrt n\, V(\gamma)\wedge
W_{n-1}(\gamma),\qquad W_{n-1}(\gamma)\in
\wedge^{n-1}(T_\gamma\Gamma_X)\end{align*} (the pairing in the
expression $\langle V(\gamma),W_n(\gamma)\rangle_\gamma$ is
carried out in the first ``variable,'' so that $a_n^*(V(\gamma))$
becomes the adjoint of $a_n(V(\gamma))$\,).

Now, for a fixed $\gamma\in\Gamma_X$, we define an operator
$R_n(\gamma)$ in $\wedge^n(T_\gamma\Gamma_X)$ as follows:
\begin{gather*} R_n(\gamma){:=}\sum_{x\in\gamma}
R(\gamma,x),\qquad
D(R_n(\gamma)){:=}\wedge_0^n(T_\gamma\Gamma_X),\\
R_n(\gamma,x){:=}\sum_{i,j,k,l=1}^d R_{i,j,k,l}(x) a^*_n(e_i)a_n
(e_j)a^*_n (e_k)a_n(e_l).\end{gather*} Here, $\{e_j\}_{j=1}^d$ is
a fixed orthonormal  basis in the space $T_xX$ considered as a
subspace of $T_\gamma\Gamma_X$, $\wedge^n_0(T_\gamma\Gamma_X)$
consists of all $W_n(\gamma)\in \wedge^n(T_\gamma\Gamma_X)$ having
only a finite number of nonzero coordinates in the direct sum
expansion \eqref{tang-n1}, and $R_{ijkl}$ is the curvature tensor
on $X$.

Next, let $A(\gamma)\in(T_{\gamma,\infty}\Gamma_X)^{\otimes 2}$,
so that $A(\gamma)=(A(\gamma,x,y))_{x,y\in\gamma}$, where
$A(\gamma,x,y)\in T_yX\otimes T_xX$. We realize $A(\gamma)$ as a
linear operator acting from $T_{\gamma,0}\Gamma_X$ into
$T_{\gamma,\infty}\Gamma_X$ by setting $$ T_{\gamma,0}\Gamma_X\ni
V(\gamma)=(V(\gamma,x))_{x\in\gamma}\mapsto
A(\gamma)V(\gamma){:=}\left(\sum_{x\in\gamma}\la
A(\gamma,x,y),V(\gamma,x)\ra_x \right)_{y\in\gamma}\in
T_{\gamma,\infty}\Gamma_X.$$ If we additionally suppose that, for
any $\Lambda\in{\cal O}_c(X)$,
$$\sum_{y\in\gamma}\left(\sum_{x\in\gamma_\Lambda}\|A(\gamma,x,y)\|_{T_yX\otimes
T_xX }\right)^2<\infty,$$ then, as easily seen, $A(\gamma)$ is
indeed an operator acting from $T_{\gamma,0}\Gamma_X$ into
$T_\gamma\Gamma_X$. In the latter case, we define a linear
operator $A(\gamma)^{\wedge n}$ in $\wedge^n(T_\gamma \Gamma_X)$
with domain $D(A(\gamma)^{\wedge
n}){:=}\wedge_0^n(T_\gamma\Gamma_X)$ as follows:
$$A(\gamma)^{\wedge n}{:=} A (\gamma)\otimes{\bf 1}\dots\otimes
{\bf 1}+ {\bf 1}\otimes
 A(\gamma) \otimes {\bf 1}\otimes\dots\otimes{\bf 1%
}+ \dots+{\bf 1}\otimes\dots\otimes{\bf 1}\otimes
 A(\gamma).
$$

We set $$B_\mu'(\gamma)=(B'_\mu(\gamma,x,y))_{x,y\in\gamma}\in
(T_{\gamma,\infty}\Gamma_X)^{\otimes 2},\qquad
B'_\mu(\gamma,x,y){:=}\nabla^X_y B_\mu(\gamma,x).$$ It follows
from \eqref{awrzejirt} that, for $\mu$-a.e.\ $\gamma\in\Gamma_X$,
\begin{align*} &
\sum_{y\in\gamma}\left(\sum_{x\in\gamma_\Lambda}\|B_\mu'(\gamma,x,y)\|_{T_yX\otimes
T_xX }\right)^2\\ &\qquad \le
\left(\sum_{y\in\gamma}\sum_{x\in\gamma_\Lambda}\|B_\mu'(\gamma,x,y)\|_{T_yX\otimes
T_xX }\right)^2<\infty.\end{align*}  Therefore, the operator
$B_\mu'(\gamma)^{\wedge n}:\wedge^n_0(T_\gamma\Gamma_X)\to
\wedge^n(T_\gamma\Gamma_X)$ is well-defined for $\mu$-a.e.\
$\gamma\in\Gamma_X$.

\setcounter{theorem}{16}
\begin{theorem}\label{awqawaewaewaeqat4w} Let the conditions of Theorem~\rom{\ref{dtdrtdtdsrt}} be
satisfied\rom. Then\rom, we have on ${\cal D}\Omega^n$\rom:
 $${\bf H}^{\mathrm
 R}_{\mu,n}W(\gamma)={\bf H}^{\mathrm B}_{\mu,n}+R_n(\gamma)W(\gamma)-B_\mu'(\gamma)^{\wedge
 n}W(\gamma),\qquad \text{\rom{$\mu$-a.e.\ }}\gamma\in\Gamma_X.$$\end{theorem}

 \noindent{\it Proof}. We fix any $W^{(1)}\in{\cal D}\Omega^n$.
  By Corollary~\ref{pouoq213w}, we have
  \begin{align}\label{12548677} {\bf H}^{\mathrm R}_{\mu,n}W^{(1)}(\gamma)&=\sum_{x\in\gamma}
({\pmb \delta}_{\mu,x,n}{\bf d}_{x,n}+{\bf
d}_{x,n}{\pmb\delta}_{\mu,x,n})W^{(1)}(\gamma)\\&\quad+
\sum_{x,y\in\gamma,\, x\ne y}({\pmb \delta}_{\mu,x,n}{\bf
d}_{y,n}+{\bf
d}_{y,n}{\pmb\delta}_{\mu,x,n})W^{(1)}(\gamma).\notag\end{align}
By \eqref{qwgfhjzjg} and \eqref{3w43276756}, we get for any
$W^{(2)}\in{\cal D}\Omega^n$
\begin{gather}\label{essresews}
\int_{\Gamma_X}\left\langle\sum_{x\in\gamma}({\pmb
\delta}_{\mu,x,n}{\bf d}_{x,n}+{\bf
d}_{x,n}{\pmb\delta}_{\mu,x,n})W^{(1)}(\gamma),W^{(2)}(\gamma)\right
\rangle_{\wedge^n(T_\gamma\Gamma_X)}\,\mu(d\gamma)\\=
\int_{\Gamma_X}\sum_{x\in\gamma}\left\langle( {\pmb
\delta}_{\mu,x,n}{\bf d}_{x,n}+{\bf
d}_{x,n}{\pmb\delta}_{\mu,x,n})W^{(1)}(\gamma),W^{(2)}(\gamma)
\right \rangle_{\wedge^n(T_\gamma\Gamma_X)}\,\mu(d\gamma)\notag\\
=\int_{\Gamma_X}\mu(d\gamma)\int_{X}\sigma(\gamma,dx)\,\langle
({\pmb \delta}_{\mu,x,n}{\bf d}_{x,n}+{\bf
d}_{x,n}{\pmb\delta}_{\mu,x,n})W^{(1)}(\gamma+\eps_x),W^{(2)}(\gamma+\eps_x)
\rangle_{\wedge^n(T_{\gamma+\eps_x}\Gamma_X)}.\notag\end{gather}
By a slight modification of the proof of the Weitzenb\"ock formula
on the manifold $X$ (see e.g.\ \cite{CFKSi}), we get for a fixed
$\gamma\in\Gamma_X$
\begin{gather}\label{6545465} \int_{X}\sigma(\gamma,dx)\,\langle
({\pmb \delta}_{\mu,x,n}{\bf d}_{x,n}+{\bf
d}_{x,n}{\pmb\delta}_{\mu,x,n})W^{(1)}(\gamma+\eps_x),W^{(2)}(\gamma+\eps_x)
\rangle_{\wedge^n(T_{\gamma+\eps_x}\Gamma_X)}\\
=\int_X\sigma(\gamma,dx)\big\langle -\Delta^X_x
W^{(1)}(\gamma+\eps_x)-\langle \nabla^X_x
W^{(1)}(\gamma+\eps_x),\beta_\sigma(\gamma,x)\rangle_x\notag\\
\text{} +R_n(\gamma+\eps_x,x)W^{(1)}(\gamma+\eps_x)-
(\nabla^X_x\beta_\sigma(\gamma,x))^{\wedge
n}W^{(1)}(\gamma+\eps_x),W^{(2)}(\gamma+\eps_x)\big\rangle_{\wedge
^n(T_{\gamma+\eps_x}\Gamma_X) }.\notag\end{gather} We note that
the function under the sign of integral on the right hand side of
equality~\eqref{6545465}, considered as a function of $\gamma$ and
$x$, is integrable with respect to the measure
$\mu^{(1)}(d\gamma,dx)$. Indeed, the integrability of the function
$$ F_1(\gamma,x){:=} \big\langle -\Delta^X_x
W^{(1)}(\gamma+\eps_x)-\langle \nabla^X_x
W^{(1)}(\gamma+\eps_x),\beta_\sigma(\gamma,x)\rangle_x,
W^{(2)}(\gamma+\eps_x)\big\rangle_{\wedge
^n(T_{\gamma+\eps_x}\Gamma_X) } $$ follows from the proof of
Theorem~\ref{th4.1}, the integrability of the function
$$F_2(\gamma,x){:=} -\big\langle
(\nabla^X_x\beta_\sigma(\gamma,x))^{\wedge
n}W^{(1)}(\gamma+\eps_x),W^{(2)}(\gamma+\eps_x)\big\rangle_{\wedge
^n(T_{\gamma+\eps_x}\Gamma_X) } $$ follows from the proof of
Theorem~\ref{dtdrtdtdsrt}, and the integrability of the function
$$F_3(\gamma,x){:=}\langle
R_n(\gamma+\eps_x,x)W^{(1)}(\gamma+\eps_x),W^{(2)}(\gamma+\eps_x)\rangle_{\wedge^n(T_{
\gamma+\eps_x}\Gamma_X)}$$ follows from the estimate $$
|F_3(\gamma,x)|\le n^2 d^4 R_\Lambda
\|W^{(1)}(\gamma+\eps_x)\|_{\wedge^n(T_{\gamma+\eps_x}\Gamma_X)}
\|W^{(2)}(\gamma+\eps_x)\|_{\wedge^n(T_{\gamma+\eps_x}\Gamma_X)}
,$$ where
$$R_{\Lambda}{:=}\max_{i,j,k,l=1,\dots,d}\sup_{x\in\Lambda}|R_{i,j,k,l}(x)|,$$
$\Lambda\in{\cal O}_c(X)$ being such that, for some compact
$\Lambda'\subset\Lambda$, $W^{(1)}(\gamma)=W^
{(1)}(\gamma_{\Lambda'})$ for all $\gamma\in\Gamma_X$.

Hence, by \eqref{qwgfhjzjg}, \eqref{essresews}, \eqref{6545465},
and Theorem~\ref{th4.1} \begin{align} \label{7862343}&
\int_{\Gamma_X}\left\langle\sum_{x\in\gamma}({\pmb
\delta}_{\mu,x,n}{\bf d}_{x,n}+{\bf
d}_{x,n}{\pmb\delta}_{\mu,x,n})W^{(1)}(\gamma),W^{(2)}(\gamma)\right
\rangle_{\wedge^n(T_\gamma\Gamma_X)}\,\mu(d\gamma)\\& \qquad =
\int_{\Gamma_X}\left\langle {\bf H}^{\mathrm
B}_{\mu,n}W^{(1)}(\gamma)+\sum_{x\in\gamma}R_n(\gamma,x)
W^{(1)}(\gamma)\right.\notag \\&\qquad\quad \left. \text{}
-\sum_{x\in\gamma}(\nabla^X_x B_\mu(\gamma,x))^{\wedge
n}W^{(1)}(\gamma),W^{(2)}(\gamma)\right\rangle_{\wedge^n(T_\gamma\Gamma_X)}\,
\mu(d\gamma) .\notag\end{align}

Next, using formulas \eqref{aweese}--\eqref{oipdrt}, we have
\begin{equation}\label{zfr352}({\pmb \delta}_{\mu,x,n}{\bf
d}_{y,n}+{\bf
d}_{y,n}{\pmb\delta}_{\mu,x,n})W^{(1)}(\gamma)=-(\nabla^X_y
B_\mu(\gamma,x))^{\wedge n}W^{(1)}(\gamma),\qquad
\gamma\in\Gamma_X,\ x,y\in\gamma,\ x\ne y.\end{equation} Thus, by
\eqref{12548677}, \eqref{7862343}, and \eqref{zfr352}, we get, for
$\mu$-a.e.\ $\gamma\in\Gamma_X$,
\begin{align*} {\bf H}^{\mathrm R}_{\mu,n}W^{(1)}(\gamma)&={\bf
H}_{\mu,n}^{\mathrm
B}W^{(1)}(\gamma,x)+R_n(\gamma)W^{(1)}(\gamma)-\sum_{x\in\gamma}
(\nabla^X_x B_\mu(\gamma,x))^{\wedge n}W^{(1)}(\gamma)\\ &\quad
-\sum_{x,y\in\gamma,\, x\ne y} (\nabla^X_y
B_\mu(\gamma,x))^{\wedge n}W^{(1)}(\gamma)\\ &={\bf
H}_{\mu,n}^{\mathrm B}W^{(1)}(\gamma)+
R_n(\gamma)W^{(1)}(\gamma)-(B_\mu'(\gamma))^{\wedge n}
W^{(1)}(\gamma).\quad\blacksquare\end{align*}



\section{Examples}
\label{dseersre}

In this section, we will discuss some measures on the
configuration space $\Gamma_X$ to which the above results are
applicable.

\subsection{(Mixed) Poisson measures}
 Let
$\pi_z$, $z>0$, denote the Poisson measure on $(\Gamma_X,{\cal
B}(\Gamma_X))$ with intensity measure $zm$.  This measure can be
characterized by its Laplace transform $$\int_{\Gamma} \exp[\la
f,\gamma\ra]\,\pi_z(d\gamma)
=\exp\bigg(\int_X(e^{f(x)}-1)\,zm(dx)\bigg),\qquad f\in\D.$$ We
refer to e.g.\ \cite{AKR1,VGG} for a detailed discussion of the
construction of the Poisson measure on the configuration space.
The measure $\pi_z$ satisfies \eqref{qwgfhjzjg} with
$\sigma(\gamma,dx)=zm(dx)$, which is the so-called Mecke identity
\cite{Mec67}.

Every measure $\pi_z$ is concentrated on the subset $\Xi_z\in{\cal
B}(\Gamma_X)$ consisting of those $\gamma\in\Gamma_X$ for which $$
\lim_{n\to\infty} \frac{|\gamma_{\Lambda_n}|}{m(\Lambda_n)}=z,$$
where $(\Lambda_n)_{n=1}^\infty$ is an extending sequence of sets
from ${\cal O}_c(X)$ such that $\Lambda_n\to X$ as $n\to\infty$,
see \cite{GGPS,NZ2}.

Let $\theta$ be a probability measure on $(0,\infty)$. A mixed
Poisson measure $\pi_\theta$ is defined by
$$\pi_\theta(\cdot){:=}\int_0^\infty \theta(dz) \pi_z(\cdot).$$
Then, evidently $\pi_\theta$ satisfies \eqref{qwgfhjzjg} with $$
\rho(\gamma,x)=z m(dx)\quad \text{for }\gamma\in\Xi_z.$$

Let us suppose that $$ \int_{0}^\infty z^n\,
\theta(dz)<\infty\qquad\text{for all }n\in\N.$$  Then, condition
\eqref{edrdt} is fulfilled, and furthermore all the theorems of
Section~\ref{dodf} are applicable to the measure $\pi_\theta$.

Let us remark the following interesting fact. The Dirichlet form
on functions, $({\cal E}_{\pi_\theta},D({\cal E}_{\pi_\theta}))$,
is irreducible if and only if $\pi_\theta$ is a pure Poisson
measure $\pi_z$, see \cite[Theorem~6.3]{AKR2}. On the other hand,
by Theorem~\ref{rtsreaweraewra}, the Bochner bilnear forms $({\cal
E}_{\pi_\theta,n}^{\mathrm B},D({\cal E}_{\pi_\theta,n}^{\mathrm
B}))$, $n\in\N$, are irreducible for all measures $\pi_\theta$. In
other words, for $\pi_\theta\ne \pi_z$ there exist
square-integrable non-constant harmonic functions, but no
square-integrable Bochner harmonic forms.

\subsection{Ruelle measures}

In this subsection, we will discuss a class of Gibbs measures on
the configuration space over $\R^d$. So, let $X{:=}\R^d$,
$d\in\N$, and let $\Gamma{:=}\Gamma_{\R^d}$. The volume measure
$m$ on $\R^d$ is now the Lebesgue measure.

A pair potential is a  measurable function $\phi\colon \R^d\to
\R\cup\{+\infty\}$ such that $\phi(-x)=\phi(x)$. We will also
suppose that $\phi(x)\in\R$ for $x\in\R^d\setminus\{0\}$. For
$\Lambda\in{\cal O}_c(\R^d)$, a conditional energy
$E_\Lambda^\phi\colon\Gamma\to {\Bbb R}\cup\{+\infty\}$ is defined
by $$
E_\Lambda^\phi(\gamma):=\begin{cases}\sum\limits_{\{x,y\}\subset\gamma,\,
\{x,y\}\cap\Lambda\ne\varnothing }\phi(x-y),&\text{if }
\sum\limits_{\{x,y\}\subset\gamma,\,
\{x,y\}\cap\Lambda\ne\varnothing}|\phi(x-y)|<\infty,\\
+\infty,&\text{otherwise}.\end{cases}$$ Given $\Lambda\in{\cal
O}_c(\R^d)$, we define for $\gamma\in\Gamma$ and $\Delta\in {\cal
B}(\Gamma)$
\begin{align*}
\Pi_\Lambda^{z,\phi}(\gamma,\Delta){:=}&{\bf
1}_{\{Z_\Lambda^{z,\phi}<\infty\}}
(\gamma)\,[Z_\Lambda^{z,\phi}(\gamma)]^{-1}\\ &\times \int_\Gamma
{\bf 1}_\Delta(\gamma_{\Lambda^c}+\gamma_\Lambda') \exp\big[
-E_\Lambda^\phi(\gamma_{\Lambda^c}+\gamma_\Lambda')
\big]\,\pi_z(d\gamma'),\end{align*} where $$
Z_\Lambda^{z,\phi}(\gamma){:=}\int_\Gamma\exp\big[
-E_\Lambda^\phi(\gamma_{\Lambda^c}+\gamma_\Lambda')
\big]\,\pi_z(d\gamma').$$

A probability measure $\mu $ on $(\Gamma,{\cal B}(\Gamma))$ is
called a grand canonical Gibbs measure with interaction potential
$\phi$ if it satisfies the Dobrushin--Lanford--Ruelle equation
$$\mu\Pi_\Lambda^{z,\phi}=\mu\qquad\text{for all }\Lambda\in{\cal
O}_c(\R^d).$$ Let ${\cal G}(z,\phi)$ denote the set of all such
probability measures $\mu$.

It can be shown \cite{Ge79} that the unique grand canonical Gibbs
measure corresponding to the free case, $\phi=0$, is the Poisson
measure $\pi_z$.

We rewrite the conditional energy $E_\Lambda^\phi$ in the
following form $$
E_\Lambda^\phi(\gamma)=E_\Lambda^\phi(\gamma_\Lambda)+W(\gamma_\Lambda\mid
\gamma_{\Lambda^c}),$$ where the term
\begin{equation}\label{111}W(\gamma_\Lambda\mid
\gamma_{\Lambda^c})=\begin{cases}\sum\limits_{x\in\gamma_\Lambda,\,y\in\gamma_
{\Lambda^c}}\phi(x-y),&\text{if
}\sum\limits_{x\in\gamma_\Lambda,\,y\in\gamma_
{\Lambda^c}}|\phi(x-y)|<\infty,\\
+\infty,&\text{otherwise},\end{cases}
\end{equation} describes the interaction energy between $\gamma_\Lambda$ and
$\gamma_{\Lambda^c}$.  Analogously, we define
$W(\gamma'\mid\gamma'')$ when $\gamma'\cap \gamma''=\varnothing$.

We suppose that the interaction potential $\phi$ is stable, that
is, the following condition is satisfied:
\begin{description}

\item[(S)] ({\it Stability})
 There exists $B\ge0$ such that, for any $\Lambda\in{\cal O}_c(\R^d)$
and for all $\gamma\in\Gamma_\Lambda$,
$$E_{\Lambda}^\phi(\gamma)\ge -B|\gamma|.$$
\end{description}

(Notice that the stability condition automatically implies that
the potential $\phi$ is semi-bounded from below.)

Then, any $\mu\in{\cal G}(z,\phi)$ satisfies  identity
\eqref{qwgfhjzjg} with
\begin{equation}\label{grtrerre}\rho(\gamma,x)=z
\exp\big[-W(\{x\}\mid\gamma)\big].\end{equation}  In fact, this
property uniquely characterizes a Gibbs measure in the sense that
any probability measure $\mu$ on $(\Gamma,{\cal B}(\Gamma))$
belongs to ${\cal G}(z,\phi)$ if and only if $\mu$ satisfies
\eqref{qwgfhjzjg} with $\rho(\gamma,x)$ given by \eqref{grtrerre}
(cf.\ \cite{NZ}, see also \cite{Kuna}).

Let us now describe a class of Gibbs measures which appears in
classical statistical mechanics of continuous systems \cite{Ru70}.
For every $r=(r^1,\dots,r^d)\in\Z^d$, we define a cube
$$Q_r:=\big\{\, x\in\R^d\mid r^i-\frac 12\le x^i<r^i+\frac12
\,\big\}.$$ These cubes form a partition of $\R^d$. For any
$\gamma\in\Gamma$, we set $\gamma_r:=\gamma_{Q_r}$, $r\in\Z^d$.
For $N\in\N$ let $\Lambda_N$ be the cube with side length $2N-1$,
centered at the origin in $\R^d$. $\Lambda_N$ is then a union of
$(2N-1)^d$ unit cubes of the form $Q_r$.

 We formulate the following conditions on
the interaction.

\begin{description}

\item[(SS)] ({\it Superstability})
 There exist $A>0$, $B\ge0$ such that if $\gamma\in\Gamma_
{\Lambda_N}$ for some $N$, then
$$E_{\Lambda_N}^\phi(\gamma)\ge\sum_{r\in\Z^d}\big[A|\gamma_r|
^2-B|\gamma_r|\big].$$

\end{description}

This condition is evidently stronger than (S).

\begin{description}

\item[(LR)] ({\it Lower regularity}) There exists a decreasing positive
function $a\colon\N\to{\Bbb R}_+$ such that
$$\sum_{r\in\Z^d}a(\|r\|)<\infty$$ and for any
$\Lambda',\Lambda''$ which are finite unions of cubes $Q_r$ and
disjoint, with $\gamma'\in\Gamma_{\Lambda'}$,
$\gamma''\in\Gamma_{\Lambda''}$,
$$W(\gamma'\mid\gamma'')\ge-\sum_{r',r''\in\Z^d}a(\|r'-r''\|)
|\gamma_{r'}'|\,|\gamma_{r''}''|.$$ Here, $\|\cdot\|$ denotes the
maximum norm on $\R^d$.

\item[(I)] ({\it Integrability}) We have
$$\int_{\R^d}|1-e^{-\phi(x)}|\,m(dx)<+\infty.$$

\end{description}

A probability measure $\mu$ on $(\Gamma,{\cal B}(\Gamma))$ is
called tempered if $\mu$ is supported by
$$S_\infty{:=}\bigcup_{n=1}^\infty S_n,$$ where $$S_n:=\big\{\,
\gamma\in\Gamma\mid \forall N\in\N\ \sum_{r\in\Lambda_N\cap\Z^d}
|\gamma_r|^2\le n^2|\Lambda_N \cap\Z^d| \,\big\}.$$ By ${\cal
G}^t(z,\phi)\subset{\cal G}(z,\phi)$ we denote the set of all
tempered grand canonical Gibbs measures (Ruelle measures for
short). Due to \cite{Ru70} the set ${\cal G}^t(z,\phi)$ is
non-empty for all $z>0$ and any potential $\phi$ satisfying
conditions (SS), (LR), and (I).

Let us now recall  the so-called Ruelle bound (cf.\ \cite{Ru70}).

\begin{theorem}\label{waessedf} Let $\phi$ be a pair potential satisfying conditions
\rom{(SS), (LR),} and \rom{(I),} and let $\mu\in{\cal
G}^t(z,\phi)$\rom, $z>0$\rom. Then\rom, for any $n\in\N$ and any
measurable symmetric function $f^{(n)}:(\R^d)^n\to[0,\infty]$, we
have
\begin{align*} &\int_\Gamma \sum_{\{x_1,\dots,x_n\}\subset\gamma}
f^{(n)}(x_1,\dots,x_n)\,\mu(d\gamma)\\&\qquad =\frac1{n!}\,
\int_{(\R^d)^n} f^{(n)}(x_1,\dots,x_n)
k_\mu^{(n)}(x_1,\dots,x_n)\,m(dx_1)\dotsm m(dx_n),\end{align*}
where $k_\mu^{(n)}$ is a non-negative measurable symmetric
function on $(\R^d)^n$\rom, called the $n$-th correlation function
of the measure $\mu$\rom, and this function satisfies the
following estimate \begin{equation}\label{534}\forall
(x_1,\dots,x_n)\in(\R^d)^n:\quad k_\mu^{(n)}(x_1,\dots,x_n)\le
\xi^n,\end{equation} where $\xi> 0$ is independent of $n$\rom.
\end{theorem}

The above theorem particularly implies that any Ruelle measure
$\mu$ satisfies \eqref{edrdt}.

We suppose:

\begin{description} \item[(S1)] There
 exists $r>0$ such that $$\int_{B(r)^c}|\phi(x)|\,
 m(dx)<\infty,$$ where $B(r)$ denotes the open
ball in $\R^d$ of radius $r$ centered at the
origin.\end{description}

\setcounter{lemma}{1}
\begin{lemma} \label{894356}
Let \rom{(SS), (LR), (I),} and \rom{(S1)} hold\rom. Then\rom,
$$\sum_{y\in\gamma}|\phi(x-y)|<\infty\qquad \text{\rom{for
$\mu\otimes m$-a.e.\ $(\gamma,x)\in\Gamma\times\R$}}.$$
Moreover\rom, for $\mu\otimes m$-a\rom.e\rom.\
$(\gamma,x)\in\Gamma\times\R^d$
$$\rho(\gamma,x)=z\exp\left[-\sum_{y\in\gamma}
\phi(x-y)\right]>0.$$
\end{lemma}

\noindent {\it Proof}. It is enough to show that, for any
$\Lambda\in{\cal O}_c(\R)$
\begin{equation}\label{hgft}
\sum_{y\in\gamma_{(\Lambda^r)^c } } |\phi(x-y)|<\infty\qquad
\text{for $\mu\otimes m$-a.e.\ }(\gamma,x)\in\Gamma\times \Lambda,
\end{equation}
where $\Lambda^r{:=}\{y\in\R^d: d(y,\Lambda)\le r\}$,
$d(y,\Lambda)$ denoting the distance from $y$ to $\Lambda$.

By Theorem~\ref{waessedf} and (S1), \begin{align*}&\int_\Gamma
\mu(d\gamma)\int_\Lambda m(dx) \sum_{y\in
\gamma_{(\Lambda^r)^c}}|\phi(x-y)|\\ &\qquad =\int_\Lambda
m(dx)\int_\Gamma \mu(d\gamma)\int_{\R^d}\gamma(dy)\, |\phi(x-y)|
{\bf 1} _{(\Lambda^r)^c}(y)\\ &\qquad =\int_\Lambda
m(dx)\int_{\R^d}m(dy)k_\mu^{(1)}(y)|\phi(x-y)| {\bf 1}
_{(\Lambda^r)^c}(y) \\ &\qquad \le \xi\int_\Lambda m(dx) \int
_{(\Lambda^r)^c} m(dy)\,|\phi(x-y)|\\ &\qquad \le\xi
m(\Lambda)\int_{B(r)^c}|\phi(y)|\, m(dy)<\infty,\end{align*} which
implies \eqref{hgft}. The second conclusion of the lemma now
trivially follows from \eqref{111} and \eqref{grtrerre}.\quad
$\blacksquare$

 We suppose also that
the two following  conditions are satisfied (compare with
\cite{AKR2}).

\begin{description} \item[(D)] ({\it Differentiability})
$e^{-\phi}$ is weakly differentiable
 on $\R^d$, $\phi$ is weakly differentiable on $\R^d\setminus\{0\}$,
 and the weak gradient $\nabla\phi$ (which is a locally $m$-integrable function on
 $\R^d\setminus
 \{0\}$) considered as an $m$-a.e.\ defined function on $\R^d$
 satisfies \begin{equation}\label{wqrte}\nabla\phi\in
 L^1(\R^d,e^{-\phi}m)\cap L^3(\R^d,e^{-\phi}m).\end{equation} \end{description}

\setcounter{remark}{2}
\begin{remark}\label{8965}\rom{ It follows from (D) that $$ \nabla
e^{-\phi}=-\nabla\phi\, e^{-\phi}\qquad \text{$m$-a.e.\ on
}\R^d.$$}\end{remark}

\begin{description} \item[(S2)] There
 exists $R>0$ such that $$\int_{B(R)^c}|\nabla\phi(x)|\,
 m(dx)<\infty.$$\end{description}

\setcounter{proposition}{3}
\begin{proposition}\label{9034} Let \rom{(SS), (LR), (I), (D), (S1)} and
\rom{(S2)} hold\rom. Then\rom, any $\mu\in{\cal G}^t(z,\phi)$\rom,
$z>0$\rom, satisfies the conditions of Theorem~\rom{\ref{th4.1}}
and
\begin{equation}\label{347634}
B_\mu(\gamma,x)=-\sum_{y\in\gamma\setminus\{x\}}\nabla\phi(x-y),\qquad
x\in\gamma,\ \text{\rom{$\mu$-a\rom.e\rom.\
}}\gamma\in\Gamma_X.\end{equation}
\end{proposition}

\noindent {\it Proof}. We first prove that, for $\mu$-a.e.\
$\gamma\in\Gamma_X$ $\rho(\gamma,\cdot)$ is weakly differentiable
on $\R^d$. We fix any $f\in{\cal D}$ and $v$, a smooth vector
field on $\R^d$ with compact support, and let $\Lambda\in{\cal
O}_c(\R^d)$ be such that the supports of both $f$ and $v$ are
contained in $\Lambda$. Let $(\Lambda_N)_{N=1}^\infty$ be the
sequence of subsets of $\R^d$ as in  (SS). Let $N\in\N$ be so big
that $\Lambda^R\subset\Lambda_N$. Then, using Remark~\ref{8965},
we get
\begin{gather}\label{rtekiuu}
 \int_{\Lambda} \exp\left[-\sum_{y\in\gamma_{\Lambda_N}}\phi(x-y)\right]
\la\nabla f(x),v(x)\ra\,zm(dx) \\= \int_\Lambda
\exp\left[-\sum_{y\in\gamma_{\Lambda_N}}\phi(x-y)\right]
f(x)\left(\sum_{y\in\gamma_{\Lambda_N}}\left\la
\nabla\phi(x-y),v(x) \right\ra
-\operatorname{div}v(x)\right)\,zm(dx)\notag\\ = \int_\Lambda
\exp\left[-\sum_{y\in\gamma_{\Lambda_N}}\phi(x-y)\right]
f(x)\left(\sum_{y\in\gamma_{\Lambda^R}}\left\la
\nabla\phi(x-y),v(x) \right\ra
-\operatorname{div}v(x)\right)\,zm(dx)\notag\\ \text{}+
\int_\Lambda
\exp\left[-\sum_{y\in\gamma_{\Lambda_N}}\phi(x-y)\right]
f(x)\left(\sum_{y\in\gamma_{\Lambda_N\setminus\Lambda^R}}\left\la
\nabla\phi(x-y),v(x) \right\ra \right)\,zm(dx). \notag\end{gather}
We know from \cite[Lemma~5.1, Proposition~5.2 and its proof]{Ru70}
that, for each $\gamma\in S_\infty$, there exists a constant
$C(\gamma)>0$ such that \begin{gather} \forall N\in\N,\ \forall
x\in\Lambda:\quad \exp\big[-W(\{x\}\mid\gamma_{\Lambda_N})\big]\le
C(\gamma).\label{eweee}\end{gather}

Moreover, analogously to the proof of \eqref{hgft}, we conclude
from (S2) that
\begin{equation}\label{wee998}\int_\Lambda
\sum_{y\in\gamma_{(\Lambda^R)^c}}|\nabla\phi(x-y)|\,zm(dx)<\infty\qquad
\text{for $\mu$-a.e.\ }\gamma\in\Gamma_X.\end{equation} Now, by
virtue of Lemma~\ref{894356}, \eqref{wqrte},
\eqref{rtekiuu}--\eqref{wee998}, and the majorized convergence
theorem, we  get
\begin{multline*} \int_{\Lambda} \exp\left[-\sum_{y\in\gamma}\phi(x-y)\right]
\la\nabla f(x),v(x)\ra\,zm(dx)\\  = \int_\Lambda
\exp\left[-\sum_{y\in\gamma}\phi(x-y)\right]
f(x)\left(\sum_{y\in\gamma}\left\la \nabla\phi(x-y),v(x) \right\ra
-\operatorname{div}v(x)\right)\,zm(dx).\end{multline*}

Therefore, for $\mu$-a.e.\ $\gamma\in\Gamma$, $\rho(\gamma,\cdot)$
is weakly differentiable on $\R^d$ and
$$\beta_\sigma(\gamma,x)=-\sum_{y \in\gamma}\nabla\phi(x-y),$$ so
that $B_\mu$ is given by \eqref{347634}.

Finally, let us show that, for any $\Lambda\in{\cal O}_c(\R^d)$,
\begin{align}\label{347645}&\int_\Gamma\left
(\sum_{x\in\gamma_\Lambda}\sum_{y\in\gamma\setminus\{x\}}|
\nabla\phi(x-y)|\right)^3\,\mu(d\gamma)\\ &\qquad=\frac18\,
\int_\Gamma\left(
\sum_{\{x,y\}\subset\gamma}|\nabla\phi(x-y)|({\bf
1}_\Lambda(x)+{\bf 1}_{\Lambda}(y))\right)^3\,
\mu(d\gamma)<\infty,\notag\end{align} which implies
\eqref{awrzwerjidrt} with $\eps=1$.

The proof of \eqref{347645} is essentially analogous to that of
\cite[Lemma~4.1]{AKR2}, so we only sketch it.  By using
\cite[Proposition~3.11]{Kuna} and Theorem~\ref{waessedf}, we get,
for any non-negative symmetric function $\varphi^{(2)}(x,y)$ on
$(\R^d)^2$: \begin{gather}\label{wqwqqww} \int_\Gamma
\left(\sum_{\{x,y\}\subset\gamma}\varphi^{(2)}(x,y)\right)^3\,\mu(d\gamma)\\
=c_1
\int_{(\R^d)^6}\varphi^{(2)}(x_1,x_2)\varphi^{(2)}(x_3,x_4)\varphi^{(2)}(x_5,x_6)
k_\mu^{(6)}(x_1,\dots,x_6)\, m(dx_1)\dotsm m(dx_6)\notag\\
\text{}+ c_2
\int_{(\R^d)^5}\varphi^{(2)}(x_1,x_2)\varphi^{(2)}(x_1,x_3)\varphi^{(2)}(x_4,x_5)
k_\mu^{(5)}(x_1,\dots,x_5)\, m(dx_1)\dotsm m(dx_5)\notag\\
\text{}+\int_{(\R^d)^4}\big(c_3\varphi^{(2)}(x_1,x_2)^2\varphi^{(2)}(x_3,x_4)
+c_4
\varphi^{(2)}(x_1,x_2)\varphi^{(2)}(x_2,x_3)\varphi^{(2)}(x_3,x_4)\notag\\
\text{} +c_5
\varphi^{(2)}(x_1,x_2)\varphi^{(2)}(x_1,x_3)\varphi^{(2)}(x_1,x_4)
\big)k_\mu^{(4)}(x_1,\dots,x_4)\, m(dx_1)\dotsm
m(dx_4)\notag\\\text{}+c_6\int_{(\R^d)^3}\big(\varphi^{(2)}(x_1,x_2)^2\varphi^{(2)}(x_1,x_3)+
c_7
\varphi^{(2)}(x_1,x_2)\varphi^{(2)}(x_1,x_3)\varphi^{(2)}(x_2,x_3)\big)\notag\\
\times k_\mu^{(3)}(x_1,x_2,x_3)\, m(dx_1)m(dx_2) m(dx_3)
\notag\\\text{} +
c_8\int_{(\R^d)^2}\varphi^{(2)}(x_1,x_2)^3k_\mu^{(2)}(x_1,x_2)\,
m(dx_1) m(dx_2),\notag\end{gather} where $c_1,\dots,c_8>0$. We
recall also the estimate (cf.\ \cite[formula~(4.29)]{AKR2})
\begin{equation}\label{rtrtrtrtt}
\forall (x_1,\dots,x_n)\in(\R^d)^n:\quad
k_\mu^{(n)}(x_1,\dots,x_n)\le R_n\exp\left[-\sum_{1\le i<j\le n
}\phi(x_i-x_j)\right],
\end{equation}
where $n\in \N$ and $R_n>0$. Finally, one proves \eqref{347645} by
using \eqref{wqrte}, \eqref{wqwqqww}, \eqref{rtrtrtrtt}, and the
semi-boundedness of the potential $\phi$ from
below.\quad$\blacksquare$

\setcounter{proposition}{4}
\begin{proposition}\label{dsresturzuok}
Let the conditions of Proposition~\rom{\ref{9034}} be
satisfied\rom, let for some ${\cal R}>0$
\begin{equation}\label{zrsews} \phi(x)\le 0,\qquad x\in
B({\cal R})^c,\end{equation} and let one of the two following
conditions is satisfied\rom:

\begin{description}

\item[\rom{(a)}] $\phi\in C(\R^d)$ and for each $\gamma\in S_\infty$
the series $\sum_{x\in\gamma}\phi(\cdot-x)$ converges locally
uniformly on $X$\rom;

\item[\rom{(b)}] $d\ge2$\rom, $\phi\in C(\R^d\setminus\{0\})$\rom,
 and for each $\gamma\in S_\infty$
 the series $\sum_{x\in\gamma}\phi(\cdot-x)$ converges locally
 uniformly on $X\setminus\gamma$\rom.

 \end{description}

Then\rom, the conditions of Theorem~\rom{\ref{rtsreaweraewra}} are
satisfied for each $\mu\in{\cal G}^t(z,\phi)$\rom.
\end{proposition}

\noindent {\it Proof}. Evidently, (a) implies condition (i) of
Theorem~\rom{\ref{rtsreaweraewra}} and (b) does (ii), so that we
only have to show \eqref{rsresare}. Let us fix any $\gamma\in
S_\infty$. It follows from the definition of $S_\infty$ that there
exists $C=C(\gamma)\in\N$ such that
\begin{equation}\label{q24}|\gamma_{\Lambda_N}|\le C
m(\Lambda_N),\qquad N\in\N.\end{equation}

Let us assume  that in \eqref{zrsews} ${\cal R}=1/4$, otherwise
only a trivial modification of the proof is needed.

For $a>0$, let $[a]$ denote the integer part of $a$. Supposing
that there exist $[\frac12\, (2N-1)^d]+1$ $Q_r$ cubes  in
$\Lambda_N$ which contain at least $3C$ points of $\gamma$, we
come to a contradiction with \eqref{q24}. Therefore, there exist
at least $(2N-1)^d-[\frac12\, (2N-1)^d]$ cubes which contain less
than $3C$ points of $\gamma$. Setting $N\to\infty$, we conclude
that there exists an infinite sequence $\{Q_{r(k)},\, k\in\N\}$ of
cubes which contain $<3C$ points of $\gamma$. Let $x_{k}$ denote
the center of the cube $Q_{r(k)}$.
 Then,
\begin{equation}\label{sewaew9876}\forall x\in B(x_k,1/4),\ k\in\N:\qquad |B(x,1/4)\cap
\gamma|< 3C.\end{equation} In case of (a), we get by
\eqref{sewaew9876}: \begin{equation}\label{ftds767665}\forall x\in
B(x_k,1/4),\ k\in\N:\qquad \sum_{y\in\gamma}\phi(x-y)\le
\operatorname{const},\end{equation} and hence
\begin{equation}\label{aawe343543}\forall x\in B(x_k,1/4),\ k\in\N:\qquad
\rho(\gamma,x)\ge \exp(-\operatorname{const}).\end{equation}
Therefore, $\sigma(\gamma,\cdot)$, as well as all measures
$\sigma^{(k)}(\gamma,\cdot)$, $k\ge2$, are infinite measures.

In the case of (b), we proceed as follows. Any ball $B(x_k,1/4)$
contains $3C$ open disjoint  balls of of radius $1/(12C)$, and at
least one of these balls does not contain any point of $\gamma$.
Therefore, each $B(x_k,1/4)$ contains a ball $B(y_k,1/(24C))$ such
that
\begin{equation}\label{daweaw}\forall x\in
B(y_k,1/(24C)):\qquad \inf_{y\in\gamma}|x-y|\ge
1/(24C).\end{equation} By (b) the function $\phi$ is bounded on
$\{x\in\R^d: 1/(24C)\le|x|\le {\cal R}\}$, and therefore by
\eqref{sewaew9876} and \eqref{daweaw}, we again conclude that all
$\sigma^{(k)}(\gamma,\cdot)$, $k\in\N$ are infinite measures.\quad
$\blacksquare$

\setcounter{proposition}{5}
\begin{proposition} \label{ews89u67} Let \rom{(SS), (LR), (I),} and
\rom{(S2)} hold\rom. Furthermore\rom, let the interaction
potential $\phi$ satisfy the following conditions\rom:

\begin{description}

\item[\rom{(i)}] $\phi\in C^2(\R^d\setminus\{0\})$\rom,
$e^{-\phi}$ is continuous on $\R^d$\rom, and $e^{-\phi}\nabla\phi$
extends to a continuous vector-valued function on $\R^d$\rom;

\item[\rom{(ii)}] for each $\gamma\in S_\infty$\rom, the series $\sum_{x\in\gamma} \phi(\cdot-
x)$\rom, $\sum_{x \in\gamma}\nabla \phi(\cdot-x)$\rom, and
$\sum_{x\in\gamma}\phi''(\cdot-x)$ converge locally uniformly on
$X\setminus\gamma$\rom;

\item[\rom{(iii)}] \eqref{wqrte} holds\rom, and furthermore\rom,
\begin{equation}\label{3244}\phi''\in
 L^1(\R^d,e^{-\phi}m)\cap L^3(\R^d,e^{-\phi}m).\end{equation}

\end{description}
Then\rom, any $\mu\in{\cal G}^t(z,\phi)$\rom, $z>0$\rom, satisfies
the conditions of Theorem~\rom{\ref{dtdrtdtdsrt}.}
\end{proposition}

\noindent {\it Proof}. As easily seen, conditions (i)--(iii) of
Theorem~\ref{dtdrtdtdsrt} are now satisfied. Indeed, let us fix
any $\gamma\in S_\infty$. By condition (ii),
\begin{equation}\label{8945}
\rho(\gamma,x)=\exp\left[-\sum_{y\in\gamma}\phi(x-y)\right]>0,\qquad
x\in\R^d\setminus\gamma.\end{equation} It follows from the
definition of $S_\infty$ that, for any $y\in S_\infty$,
$\gamma\setminus\{y\}$ again belongs to $S_\infty$, and therefore,
the function $${\cal O}_{\gamma,y}\ni x\mapsto
\exp\left[-\sum_{z\in\gamma}\phi(x-z)\right]=\exp[-\phi(x-y)]
\exp\left[-\sum_{z\in\gamma\setminus\{y\}}\phi(x-z)\right]$$ is
continuous by (i) and (ii). Hence, $\rho(\gamma,\cdot)$ is
continuous on $\R^d$. Moreover, by (i), (ii), and \eqref{8945},
the function $\rho(\gamma,\cdot)$ is two times differentiable on
$\R^d\setminus\gamma$, and analogously to the above, we conclude
that the form \begin{multline*} {\cal O}_{\gamma,y} \ni x\mapsto
\nabla_x\rho(\gamma,x)=-\exp\left[-\phi(x-y)
-\sum_{z\in\gamma\setminus\{y\}}\phi(x-z)\right]\\ \times
\bigg(\nabla\phi(x-y) +\sum_{z\in
\gamma\setminus\{y\}}\nabla\phi(x-z)\bigg)
\end{multline*} is continuous on ${\cal O}_{\gamma,y}$, so that $\nabla_x\rho(\gamma,\cdot
)$ is continuous on $\R^d$. Finally, for any
$x\in\R^d\setminus\gamma$, $$\nabla_x \rho(\gamma+\eps_y,x)=
-\exp\left[-\phi(x-y)-\sum_{z\in\gamma}\phi(x-z)\right]\bigg(
\nabla\phi(x-y)+ \sum_{z\in\gamma}\nabla\phi(x-z)\bigg)$$ is
differentiable in $y$ on $\R^d\setminus(\gamma\cup\{x\})$, and $$
\frac{\rho(\gamma+\eps_x,y)}{\rho(\gamma+\eps_y,x)}\,\nabla_x\rho(\gamma+\eps_y,x)
=-\exp\left[-\phi(x-y)-\sum_{z\in\gamma}\phi(z-y)\right]\bigg(
\nabla\phi(x-y)+ \sum_{z\in\gamma}\nabla\phi(x-z)\bigg)$$ extends
to a continuous form in $y$ on $\R^d$.

That \eqref{awrzwerjidrt} holds follows from \eqref{wqrte} and
(S2) (see the proof of Proposition~\ref{9034}). Thus, it only
remains to show that \eqref{awrzejirt} is also satisfied.

It follows from the above that, for each $\gamma\in S_\infty$,
$$B_\mu(\gamma,x)=-\sum_{y\in\gamma\setminus\{x\}}\nabla\phi(x-y),\qquad
x\in\gamma, $$ and hence, by (i), (ii), we get for any
$x,y\in\gamma$: $$\nabla_y
B_\mu(\gamma,x)=\begin{cases}\phi''(x-y),&\text{if }x\ne y,\\
-\sum_{z\in\gamma\setminus\{x\}}\phi''(x-z),&\text{if
}x=y.\end{cases}$$ Hence, for any $\Lambda\in{\cal O }_c(\R^d)$,
we get
\begin{gather}\label{rerer}\int_\Gamma\left(
\sum_{y\in\gamma}\sum_{x\in\gamma_\Lambda}\|\nabla_y
B_\mu(\gamma,x)\|\right)^3\,\mu(d\gamma)\\ = \int_\Gamma\left(
\sum_{x\in\gamma_\Lambda}\|\nabla_x
B_\mu(\gamma,x)\|+\sum_{x\in\gamma_\Lambda}\sum_{y\in\gamma\setminus\{x\}}\|\nabla_y
B_\mu(\gamma,x)\|\right)^3\,\mu(d\gamma)\notag\\ \le \left( 2
\sum_{x\in\gamma_\Lambda}\sum_{y\in\gamma\setminus\{x\}}\|\phi''(x-y)\|\right)^3\,\mu(
d\gamma).\notag\end{gather} The finiteness of the latter integral
in \eqref{rerer} follows from \eqref{3244} in the same way as
\eqref{347645} follows from \eqref{wqrte}.\quad $\blacksquare$

\setcounter{remark}{6}
\begin{remark}\rom{ Let the interaction potential $\phi$ satisfy
conditions of Propositions~\ref{ews89u67}. Then, by using
Lemma~\ref{waedfdrtdtgfz}, Theorem~\ref{awqawaewaewaeqat4w}, and
Proposition~\ref{9034}, we easily see that, for every $W\in{\cal D
}\Omega^1$,
\begin{gather*}
 {\cal E}_{\mu,1}^{\mathrm R} (W,W)=
\int_\Gamma \mu(d\gamma)\int_{{\Bbb
R}^d}m(dx)\,\exp\bigg(-\sum_{y\in\gamma}\phi(x-y)\bigg)\,
|\nabla_x W(\gamma+\varepsilon_x,x)|^2\\\text{}+\frac12\,
\int\mu(d\gamma)\int_{{\Bbb R}^d}m(dx) \int_{{\Bbb R}^d}m(dy)\,
\exp\bigg(-\sum_{x'\in\gamma}\phi(x-x')-\sum_{y'\in\gamma}\phi(y-y')-\phi(x-y)
\bigg)\\ \times \phi''(x-y)
\big(W(\gamma+\varepsilon_x+\varepsilon_y,x)-W(\gamma+\varepsilon_x+\varepsilon_y,y)\big)
\big(W(\gamma+\varepsilon_x+\varepsilon_y,x)-W(\gamma+\varepsilon_x+\varepsilon_y,y)\big).
\end{gather*}
}\end{remark}

Finally, we present several examples of potentials which satisfy
conditions of Propositions~\ref{9034}--\ref{ews89u67}.

{\bf Example 1.} $\phi\in C^2_0(\R^d)$, $\phi\ge0$ on $\R^d$, and
$\phi(0)>0$.

{\bf Example 2.} ({\it Lennard--Jones type potentials}) $\phi\in
C^2(\R^d\setminus\{0\})$, $\phi\ge0$ on $\R^d$,
$\phi(x)=c|x|^{-\alpha}$ for $x\in B(r_1)$, $\phi(x)=0$ for $x\in
B(r_2)^c$, where $c>0$, $\alpha>0$, $0<r_1<r_2<\infty$.

{\bf Example 3.} ({\it Lennard--Jones \rom{6--12} potentials})
$d=3$, $\phi(x)=c(|x|^{-12}-|x|^{-6})$, $c>0$.

\subsection{Gibbs measures on configuration spaces over manifolds}

In this subsection, we will shortly discuss the case of a Gibbs
measure $\mu$ on $\Gamma_X$, where $X$ is again a general
manifold.

We formulate the following conditions on the interaction potential
$\phi$, which is now a symmetric functions
$\phi:X^2\to\R\cup\{+\infty\}$.

\newcommand{\esssup}{\operatornamewithlimits{ess\,sup}}

\begin{description}

\item[(S)] ({\it Stability})
 There exists $B\ge0$ such that, for any $\Lambda\in{\cal O}_c(X)$
and for all $\gamma\in\Gamma_\Lambda$,
$$E_{\Lambda}^\phi(\gamma){:=}\sum_{\{x,y\}\subset\gamma}\phi(x,y)\ge
-B|\gamma|.$$

\item[(I)] ({\it Integrability})
We have $$C{:=}\esssup_{x\in X}\int _X
|e^{-\phi(x,y)}-1|\,m(dy)<\infty.$$

\item[(F)] ({\it Finite range}) There exists $R>0$ such that $$
\phi(x,y)=0\quad\text{if }d(x,y)\ge R.$$

\end{description}

In a completely analogous way as for  the case of $\R^d$, one
defines a Gibbs measure $\mu$ corresponding to the interaction
potential $\phi$ and activity parameter $z>0$, and one denotes  by
${\cal G}(z,\phi)$ the set of all such measures.

\setcounter{theorem}{7}
\begin{theorem}[\cite{KKS98,KunaPhD,K}]\label{gdtrsresa}
\rom{1)} Let \rom{(S), (I),} and \rom{(F)} hold\rom, and let $z>0$
be such
 that $$z<\frac1{2e}\,( e^{2 B}C)^{-1},$$ where
$B$ and $C$ are as in \rom{(S)} and \rom{(I),} respectively\rom.
Then\rom, there exists a Gibbs measure $\mu\in{\cal G}(z,\phi)$
such that the correlation functions $k_\mu^{(n)}$ of the measure
$\mu$ satisfy the Ruelle bound \eqref{534}\rom.

\rom{2)} Let $\phi$ be a non-negative potential which fulfills
\rom{(I)} and \rom{(F).} Then\rom, for each $z>0$\rom, there
exists a Gibbs measure $\mu\in{\cal G}(z,\phi)$ such that the
correlation functions $k_\mu^{(n)}$ of the measure $\mu$ satisfy
the Ruelle bound \eqref{534}\rom.

\end{theorem}

\setcounter{proposition}{8}
\begin{proposition} \label{ewserw9854} Suppose the conditions of
Theorem~\rom{\ref{gdtrsresa}} are satisfied  and furthermore
suppose that the interaction potential $\phi$ satisfies the
following conditions\rom:

\begin{description}

\item[\rom{(i)}] $\phi\in C^2(X^2\setminus \widetilde X^2)$\rom,
$e^{-\phi}$ is continuous on $X^2$\rom, and
$e^{-\phi}\nabla^X_1\phi$ extends to a continuous vector field on
$X^2$ \rom(here $\nabla^X_1\phi$ denotes the gradient of the
function $\phi$ in the first variable\rom)\rom;

\item[\rom{(ii)}] we have $$ \esssup_{x\in X}\int_X
\big|(\nabla^X_x)^k\phi(x,y)\big|^n\exp\big(-\phi(x,y)\big)\,
m(dy)<\infty,\qquad k=1,2,\ n=1,2,3.$$

\end{description}
Let $\mu\in{\cal G}(z,\phi)$ be as in
Theorem~\rom{\ref{gdtrsresa}.} Then\rom, $\mu$ satisfies the
conditions of Theorems~\rom{\ref{th4.1}} and
\rom{\ref{dtdrtdtdsrt}.}

\end{proposition}

\noindent {\it Proof}. The proof of this proposition essentially
follows the lines of the proof of Proposition~\ref{ews89u67}, and
is even easier, since due to condition (F) all series
$\sum_{y\in\gamma}\phi(x,y)$, $\gamma\in\Gamma_X$, $x\in
X\setminus\gamma$, are finite.\quad $\blacksquare$

\setcounter{proposition}{9}
\begin{proposition}\label{awqaw98}
Suppose that the manifold $X$ satisfies the following
condition\rom:
 \begin{equation}\label{7589765} \forall r>0:\qquad 0<\inf_{x\in X}
m(B(x,r))\le \sup_{x\in X} m(B(x,r))<\infty.\end{equation} Assume
that the conditions of Propositions~\rom{\ref{ewserw9854}} are
satisfied  and either $\phi$ is a continuous bounded function on
$X^2$\rom, or $d\ge 2$ and $$ \forall r>0:\qquad \sup_{x\in
X}\sup_{y\in X,\, d(x,y)\ge r}|\phi(x,y)|<\infty.$$ Let
$\mu\in{\cal G}(z,\phi)$ be as in Theorem~\rom{\ref{gdtrsresa}.}
Then\rom, the conditions of Theorem~\rom{\ref{rtsreaweraewra}} are
satisfied\rom.
\end{proposition}

\setcounter{remark}{10}
\begin{remark}\rom{Condition \eqref{7589765}    is satisfied in the case of a manifold
having bounded geometry, see \cite{Davies}. The upper estimate
$\sup_{x\in X}m(B(x,r))<\infty$, $r>0$, holds for manifolds having
non-negative Ricci curvature (see e.g.\
\cite[Proposition~5.5.1]{Davies}). }
\end{remark}

\noindent {\it Proof}. Let us fix any sequence $\{B(x_n,2R),\,
n\in\N\}$ of disjoint balls in $X$, where  $R$ is as in (F). Let
\begin{equation}\label{3243} \Lambda_N{:=}\bigcup_{n=1}^N
B(x_n,2R),\qquad N\in\N. \end{equation} By \eqref{7589765},
$m(\Lambda_N)\to\infty$ as $N\to\infty$. By
Theorem~\ref{gdtrsresa}, the correlation functions $k_\mu^{(n)}$
satisfy the Ruelle bound. Hence, it follows from (the proof of)
\cite[Theorem~2.5.4]{Kuna}  that there exists a subsequence
$\{\Lambda_{N(k)},\, k\in\N\}$ such that, for $\mu$-a.e.\
$\gamma\in\Gamma$, there exists $C=C(\gamma)>0$ satisfying
\begin{equation}\label{4332} |\gamma_{\Lambda_{N(k)}}|\le C
m(\Lambda_{N(k)})\qquad \text{for all }k\in\N.\end{equation} By
\eqref{7589765}--\eqref{4332}, $$ |\gamma_{\Lambda_{N(k)}}|\le
 C \big(\sup_{x\in X}m(B(x,2R))\big) N(k),\qquad k\in\N.$$
Since by \eqref{7589765} $\inf_{x\in X}m(B(x,r))>0$, $r>0$,  the
rest of the proof is now completely analogous to the proof of
Proposition~\ref{dsresturzuok}. \quad $\blacksquare$

{\bf Example.} Suppose that the manifold $X$  satisfies
\eqref{7589765}, and for some $R>0$ $$\sup_{x\in X}\sup_{y\in
B(x,R)}|\nabla^X_y f(x,y)|_{T_yX}<\infty,\qquad k=1,2,$$ where
$$X^2\ni(x,y)\mapsto f(x,y){:=}d(x,y)^2\in\R.$$ (For example,
these conditions are satisfied if the manifold has a periodical
structure.) Let $\Phi\in C^2([0,\infty))$ be such that $\Phi\ge 0$
on $[0,\infty)$ and $\Phi(x)=0$ for $x\ge R^2$. Then, the
potential $\phi(x,y){:=}\Phi(f(x,y))$ satisfies the conditions of
Propositions~\ref{ewserw9854}, \ref{awqaw98}.

\vspace{5mm}

\begin{center}
 Institut f\"{u}r Angewandte Mathematik,
Universit\"{a}t Bonn, Wegelerstr.~6, D-53115 Bonn, Germany;
 SFB 256, Univ.\ Bonn, Germany; SFB 237, Bochum--D\"usseldorf--Essen, Germany;
  BiBoS, Univ.\ Bielefeld, Germany;
 CERFIM (Locarno);
Acc.\ Arch.\ (USI), Switzerland\\ e-mail:
albeverio$@$uni-bonn.de\\[3mm]Nottingham Trent University, Burton
Street, Nottingham NG1 4BU, U.K.; SFB 256, Univ. Bonn, Germany;
  BiBoS, Univ.\ Bielefeld, Germany;  Institute of Mathematics, Kiev,
  Ukraine\\ e-mail: alexei.daletskii$@$ntu.ac.uk\\[3mm] Fakult\"at f\"ur Mathematik, Universit\"at
Bielefeld, Postfach 10 01 31, D-33501 Bielefeld, Germany;
 SFB 256, Univ.\ Bonn, Germany;
 BiBoS, Univ.\ Bielefeld, Germany; Institute of
Mathematics, Kiev, Ukraine\\ e-mail:
kondrat$@$mathematik.uni-bielefeld.de\\[3mm] Institut f\"{u}r
Angewandte Mathematik, Universit\"{a}t Bonn, Wegelerstr.~6,
D-53115 Bonn, Germany;
 SFB 256, Univ.\ Bonn, Germany;
 BiBoS, Univ.\ Bielefeld, Germany\\ e-mail:
 lytvynov@wiener.iam,uni-bonn.de
\end{center}

\end{document}